\documentclass[a4paper,12pt,reqno]{amsart}
\usepackage{amsmath,graphicx,graphics,amssymb}
\newtheorem{theo}{Theorem}[section]

\newtheorem{cor}[theo]{Corollary}

\numberwithin{equation}{section}

\newcommand{\M}{\operatorname{M}}

\newcommand{\h}{\operatorname{H}}

\newcommand{\Z}{\mathbb{Z}}

\newmuskip\pFqskip
\mathchardef\pFcomma=\mathcode`, 

\begin{document}

\title{Cruciform regions and a conjecture of Di Francesco}

\author{Mihai Ciucu}
\address{Department of Mathematics, Indiana University, Bloomington, Indiana 47405}

\thanks{Research supported in part by Simons Foundation Collaboration Grant 710477}

\begin{abstract} A recent conjecture of Di Francesco states that the number of domino tilings of a certain family of regions on the square lattice is given by a product formula reminiscent of the one giving the number of alternating sign matrices. These regions, denoted ${\mathcal T}_n$, are obtained by starting with a square of side-length $2n$, cutting it in two along a diagonal by a zigzag path with step length two, and gluing to one of the resulting regions half of an Aztec diamond of order $n-1$. Inspired by the regions ${\mathcal T}_n$, we construct a family $C_{m,n}^{a,b,c,d}$ of cruciform regions generalizing the Aztec diamonds and we prove that their number of domino tilings is given by a simple product formula. Since (as it follows from our results) the number of domino tilings of ${\mathcal T}_n$ is a divisor of the number of tilings of the cruciform region $C_{2n-1,2n-1}^{n-1,n,n,n-2}$, the special case of our formula corresponding to the latter can be viewed as partial progress towards proving Di Francesco's conjecture.
\end{abstract}

\maketitle

%
%

\section{Introduction}

\begin{figure}[t]
\centerline{
\hfill
{\includegraphics[width=0.25\textwidth]{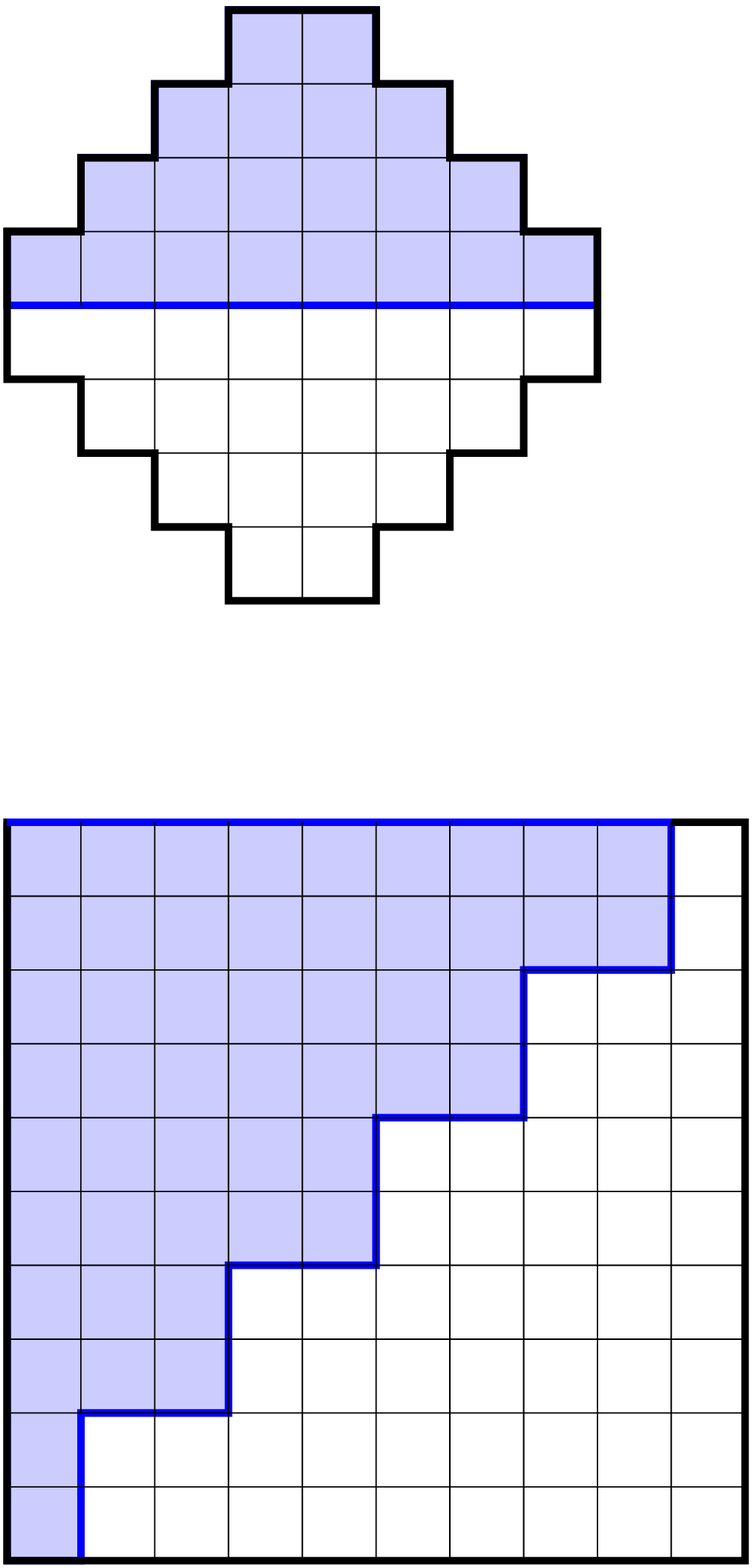}}
\hfill
{\includegraphics[width=0.22\textwidth]{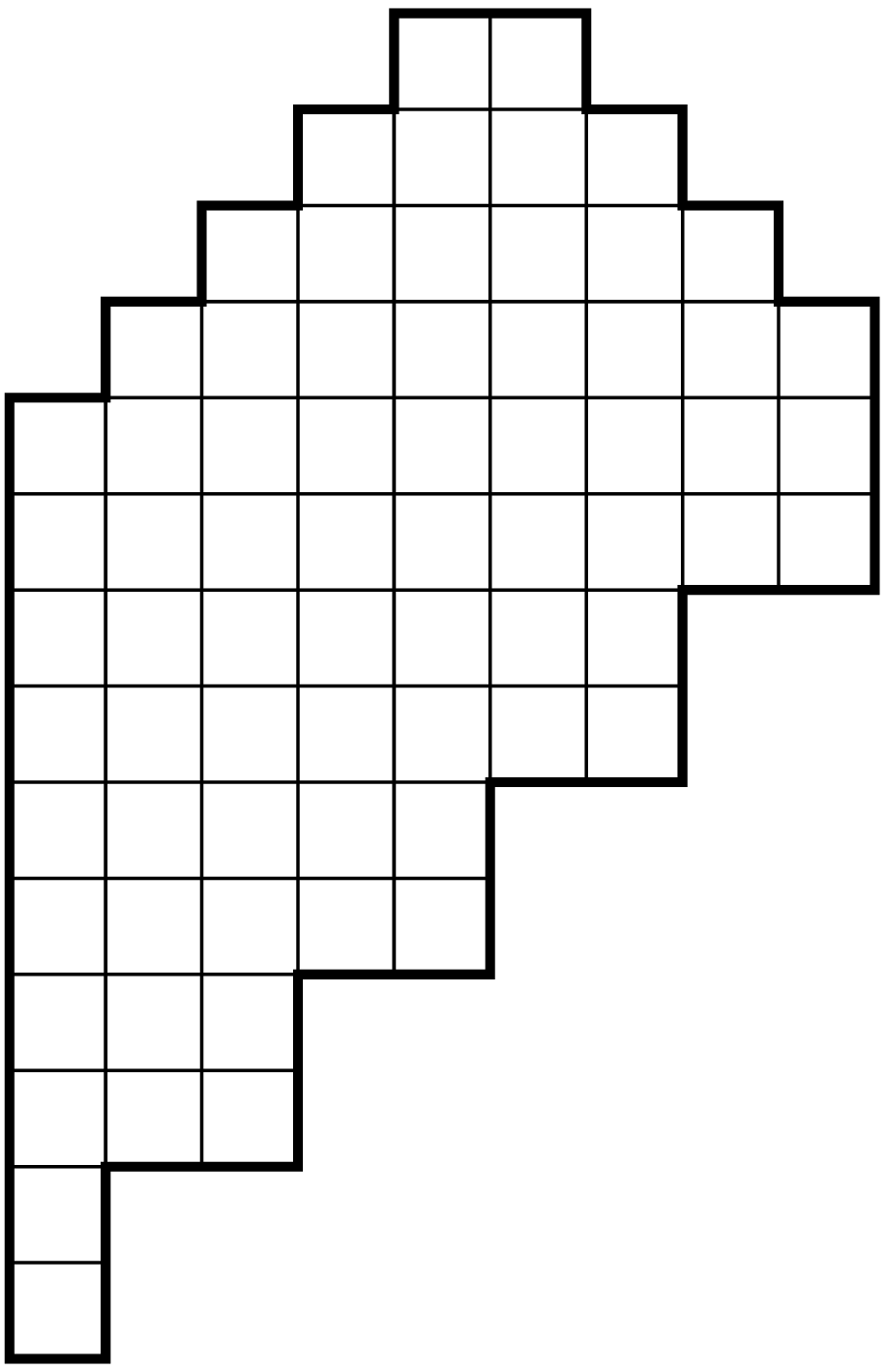}}
\hfill
}
\vskip0.1in
\caption{Half of the square of side length $2n$, and half of the Aztec diamond of order $n-1$, for $n=5$ (left). The region ${\mathcal T}_n$ for $n=5$ (right).}
\vskip-0.1in
\label{faa}
\end{figure}

In the recent paper \cite{DF} Di Francesco considers the following family of regions ${\mathcal T}_n$ on the square lattice. Let $S_{2n}$ be the lattice square of side length $2n$. Cut $S_{2n}$ into two congruent parts along its diagonal parallel to the first bisector by a zigzag line with steps of length two, leaving the bottom left unit square above the cut; denote the region above the cut by $HS_{2n}$. Let $AD_{n-1}$ be the Aztec diamond region of order $n-1$, and denote by $HD_{n-1}$ its top half. Then ${\mathcal T}_n$ is obtained by placing $HD_{n-1}$ on top of $HS_{2n}$ so that they are right-justified (see Figure \ref{faa}). Thus,  ${\mathcal T}_n$ is a genuine hybrid between a lattice square and an Aztec diamond.

Di Francesco conjectured \cite[Conjecture 8.1]{DF} that the number of domino tilings of ${\mathcal T}_n$ is equal to
\begin{equation}
2^{n(n-1)/2}\prod_{i=0}^{n-1}\frac{(4i+2)!}{(n+2i+1)!},
\label{eaa}  
\end{equation}
a formula that, as he points out, is reminiscent of the one giving the number $A_n$ of alternating sign matrices of order $n$, which is $A_n=\prod_{i=0}^{n-1}\frac{(3i+1)!}{(n+i)!}$ (see \cite{Zeil,Kup,Fischer} for three different proofs). One puzzling aspect of formula \eqref{eaa} is that the ``fingerprint'' of the Aztec diamond half of ${\mathcal T}_n$ is clearly visible in the prefactor (as the number of domino tilings of $AD_{n-1}$ is $2^{n(n-1)/2}$), but the effect of the presence of the square half $HS_{2n}$ --- whose number of domino tilings is, as a consequence of the Temperley-Fisher-Kasteleyn formula \cite{TF,Fisher-dimer,Kast} and the factorization theorem of \cite{FT}, equal to
\begin{equation}
2^{n(n-1)/2}\sqrt{\prod_{j=1}^{n}\prod_{k=1}^{n}\left(\cos^2\frac{\pi j}{2n+1}+\cos^2\frac{\pi k}{2n+1}\right)}
\label{eab}  
\end{equation}
--- is harder to understand, as the product in \eqref{eaa} that appears in its stead is much nicer (no doubt, a consequence of hybridization).

Di Francesco's interest in the regions ${\mathcal T}_n$ stems from their close connection with the twenty vertex model, a lattice model in statistical physics whose states are orientations of edges of subgraphs of the triangular lattice in which around each vertex three of the six incident edges point in and three point out (the fact that there are ${6\choose3}=20$ such choices for a given vertex gives the model its name). These connections are explained in detail in \cite{DF}.

In this paper, guided by the fact that the factorization theorem of \cite{FT} expresses the number of perfect matchings of a symmetric planar graph as the product of the number of perfect matchings of its two ``halves,'' starting with the region ${\mathcal T}_n$ we proceed to symmetrizing it by constructing a region $W_n$ that produces ${\mathcal T}_n$ as one of the halves resulting when applying the factorization theorem. It turns out that the other half also has only small primes in the factorization of its number of domino tilings, and thus --- by the factorization theorem --- so does $W_n$. We then repeat this symmetrizing process two more times, arriving at a cruciform region $C_n$ which only has simple zigzag boundaries (i.e., of step length 1).

The happy circumstance is that the number of domino tilings of $C_n$
also turns out to have only small primes in its factorization. In fact, a generalization of the regions $C_n$ involving several additional parameters, namely what we call the cruciform regions $C_{m,n}^{a,b,c,d}$, turns out to possess the same property. The main result of this paper is Theorem \ref{tba}, in which we prove a formula in the style of~\eqref{eaa} for the number of domino tilings of $C_{m,n}^{a,b,c,d}$. This is a very natural generalization of the Aztec diamond theorem of Elkies, Kuperberg, Larsen and Propp \cite{EKLP} which, surprisingly, went overlooked until now.

Having resulted from a sequence of three symmetrizations, the region $C_n$ --- which, to be specific, is the cruciform region $C_{2n-1,2n-1}^{n-1,n,n,n-2}$ --- has ${\mathcal T}_n$ as one of its eight fundamental regions. By our constructions and by the factorization theorem of \cite{FT}, it follows that the number of domino tilings of ${\mathcal T}_n$  is a divisor of our explicit product formula for the number of tilings of $C_{2n-1,2n-1}^{n-1,n,n,n-2}$. This can be viewed therefore as partial progress towards proving Di Francesco's conjecture \eqref{eaa}.

One way to prove the conjecture would be to work backwards from the cruciform region $C_{2n-1,2n-1}^{n-1,n,n,n-2}$, and identify the two factors resulting when applying the factorization theorem. Three such steps would lead us back  to the region ${\mathcal T}_n$, thus proving conjecture \eqref{eaa}. In this paper we make the first of these three steps, resulting in the enumeration of  domino tilings of the top half of the cruciform region $C_{n,n}^{a,b,b,a}$ (see Theorem \ref{tbb}).

\section{Statement of main results}

\begin{figure}[t]
\centerline{
\hfill
{\includegraphics[width=0.30\textwidth]{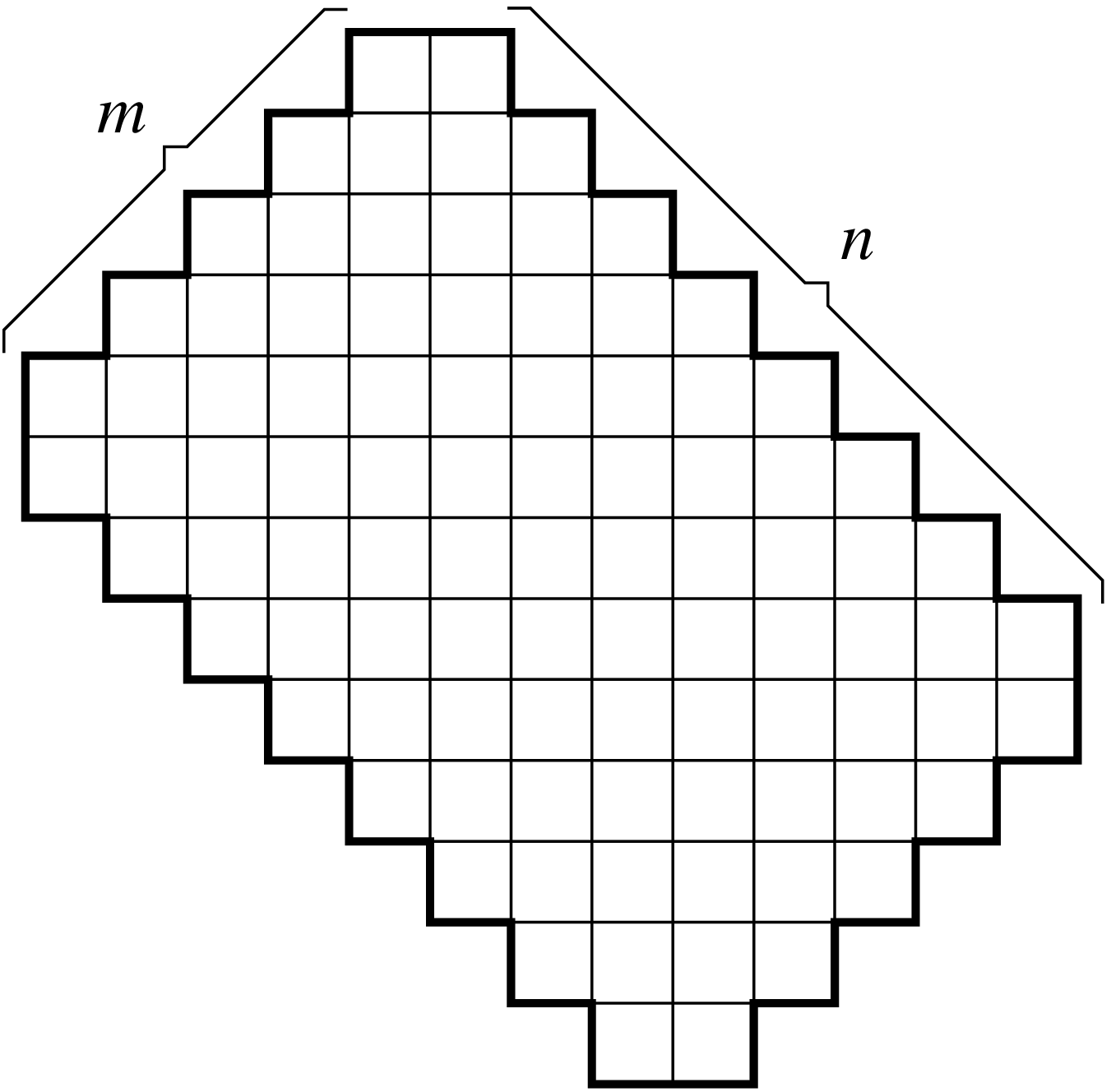}}
\hfill
{\includegraphics[width=0.20\textwidth]{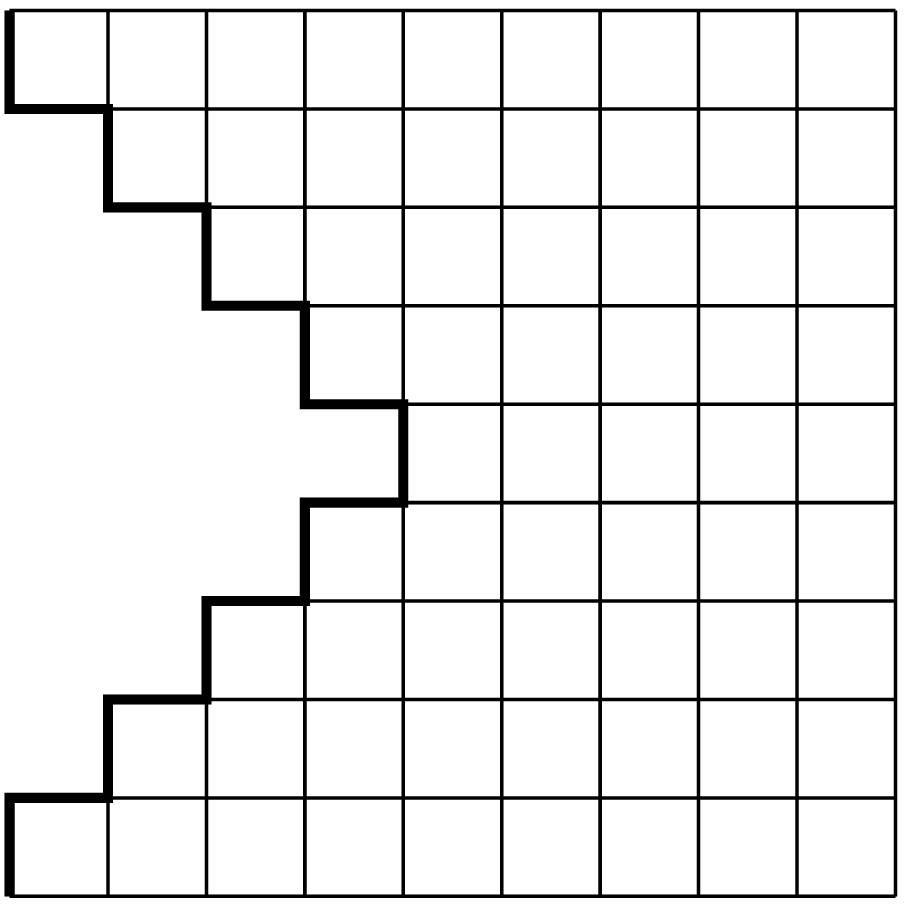}}
\hfill
}
\vskip-0.1in
\caption{The Aztec diamond region $AR_{m,n}$ for $m=5$, $n=8$ (left). Allowed corner type formed when superimposing two Aztec rectangles (right; its rotations by multiples of $90^\circ$ are also allowed).}
\vskip-0.1in
\label{fba}
\end{figure}

Recall that for $m,n\geq1$, the Aztec rectangle region $AR_{m,n}$ is the region shown on the left in Figure \ref{fba} (when $m=n$, it becomes the Aztec diamond $AD_n$). One way to define it\footnote{ This elegant definition is due to Ken Fan.} is to consider a $(2m+1)\times(2n+1)$ lattice rectangle, color its unit squares black and white in a chessboard fashion so that the corners are black, and define $G_{m,n}$ to be the graph whose vertices are the white squares, with an edge connecting two vertices precisely if the corresponding unit squares are diagonally adjacent. Then $AR_{m,n}$ is the lattice region whose planar dual is the graph $G_{m,n}$. This definition has the advantage that it shows what $AR_{m,n}$ is if $m$ and $n$ are allowed to equal zero. In particular, $AR_{0,n}$ consists of a string of $n$ diagonally adjacent unit squares.

\begin{figure}[t]
\centerline{
\hfill
{\includegraphics[width=0.40\textwidth]{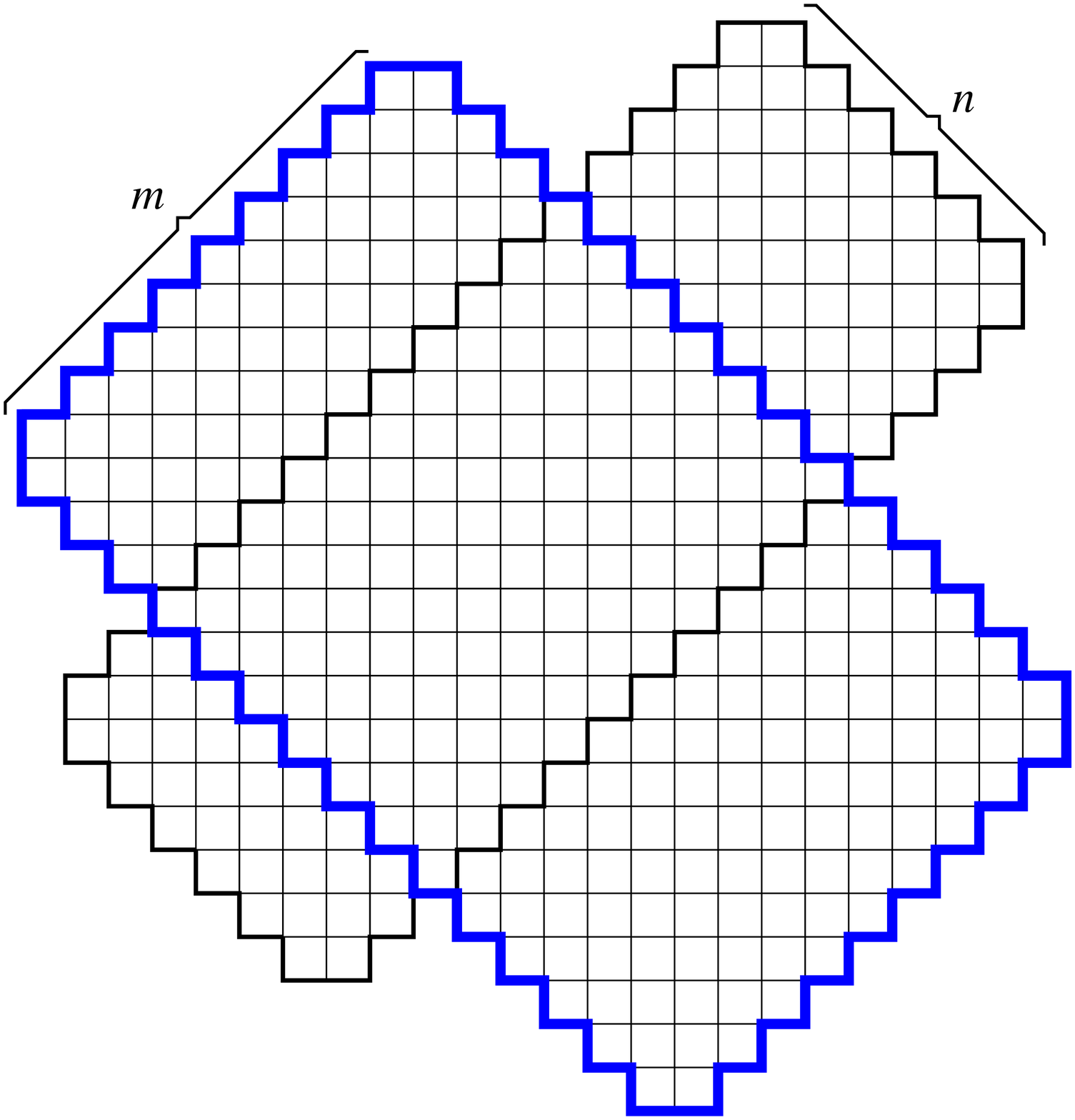}}
\hfill
{\includegraphics[width=0.40\textwidth]{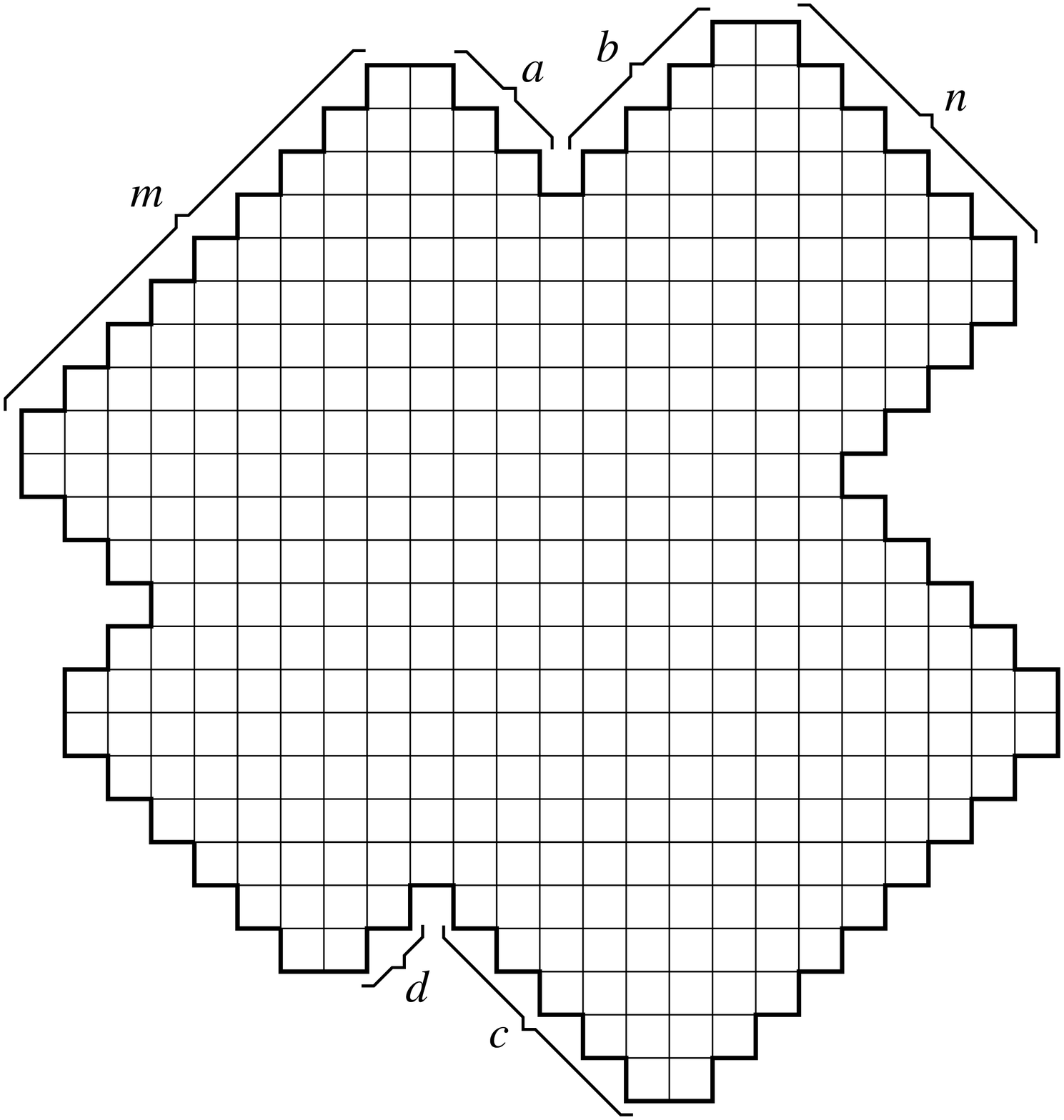}}
\hfill
}
\vskip-0.1in
\caption{Superimposing two Aztec rectangle regions so that only outer corners of the allowed type are formed (left). The cruciform region $C_{m,n}^{a,b,c,d}$ for $m=9$, $n=6$, $a=3$, $b=4$, $c=5$, $d=2$ (right).}
\vskip-0.1in
\label{fbb}
\end{figure}

To define our cruciform regions, consider two Aztec diamond regions, $AR_{m,x}$ and $AR_{y,n}$, and superimpose them in such a way that four outer corners of the type shown on the right in Figure~\ref{fba} are formed (see the picture on the left in Figure \ref{fbb} for such an example of superimposition). If in the resulting  cruciform region the four ``piers'' on the northwest, northeast, southeast and southwest stick out $a$, $b$, $c$ and $d$ units, respectively (in the example on the right in Figure \ref{fbb} we have $a=3$, $b=4$, $c=5$ and $d=2$), we denote it by $C_{m,n}^{a,b,c,d}$. This defines the cruciform regions for $a,b,c,d\geq0$.

\begin{figure}[t]
\centerline{
\hfill
{\includegraphics[width=0.44\textwidth]{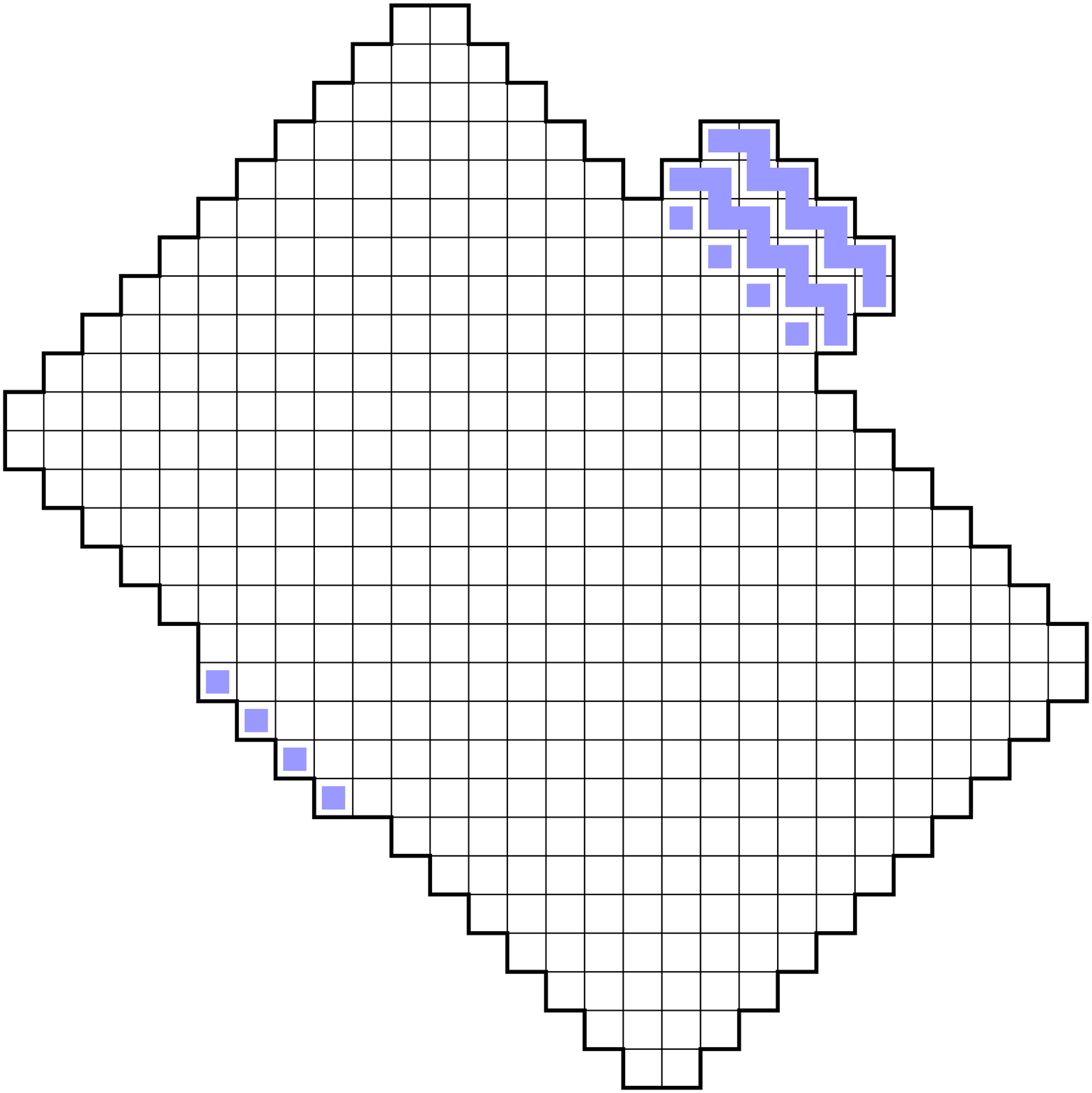}}
\hfill
{\includegraphics[width=0.53\textwidth]{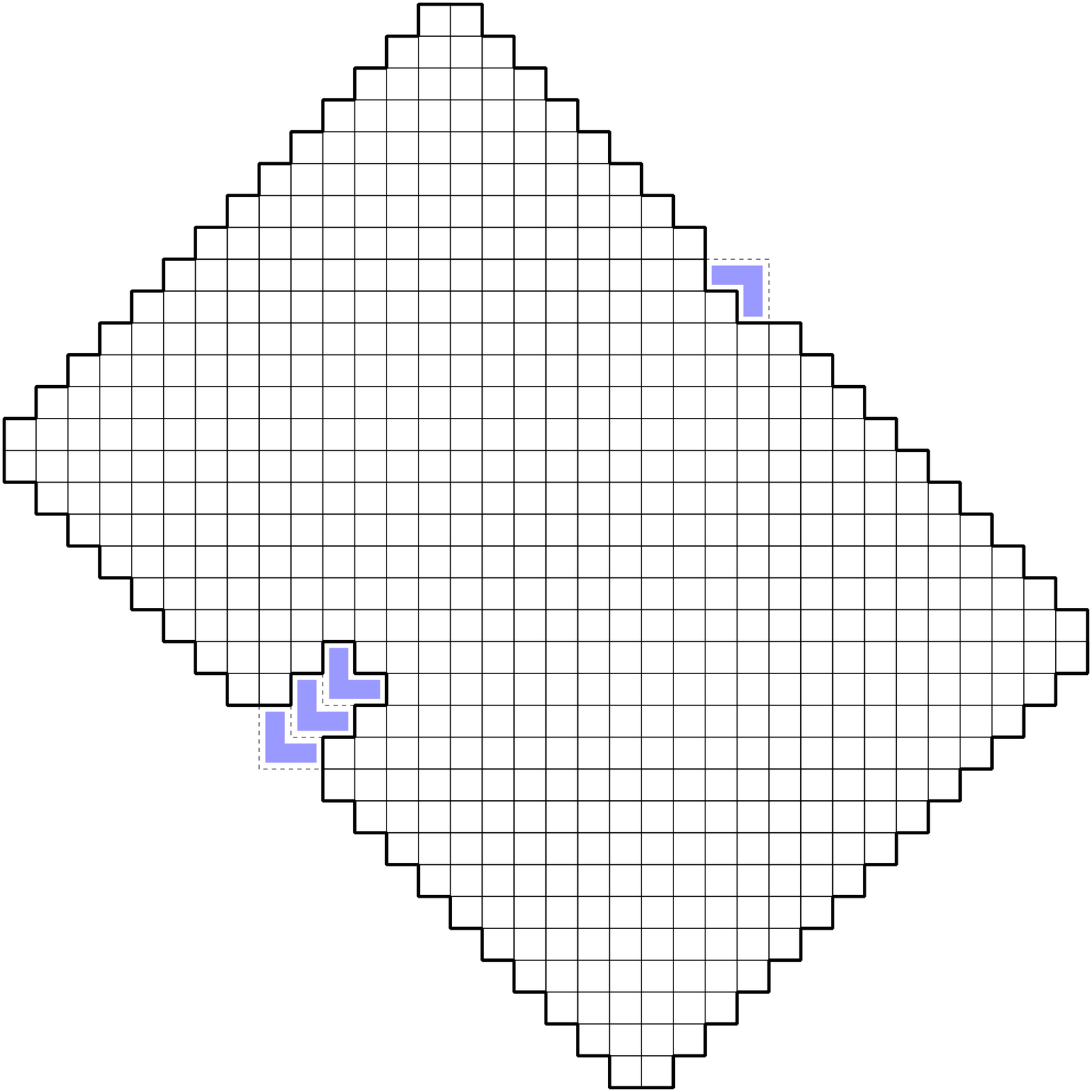}}
\hfill
}
\vskip-0.1in
\caption{{\it Left:} The cruciform region $C_{11,4}^{5,2,7,0}$; the southwestern pier, corresponding to $d=0$, consists just of 4 diagonally adjacent unit squares; the northeastern pier (corresponding to $b=2$) can be viewed as being obtained by extending the $b=0$ pier with 2 zigzag strips of unit width. {\it Right:} The cruciform region $C_{14,1}^{8,-1,10,-3}$; the ``piers'' of negative length $b=-1$ and $d=-3$ are actually bays, obtained by removing $|b|$, resp $|c|$ zigzag strips from the $b=0$, resp. $c=0$ piers.}
\vskip-0.1in
\label{fbc}
\end{figure}

\begin{figure}[t]
\centerline{
\hfill
{\includegraphics[width=0.44\textwidth]{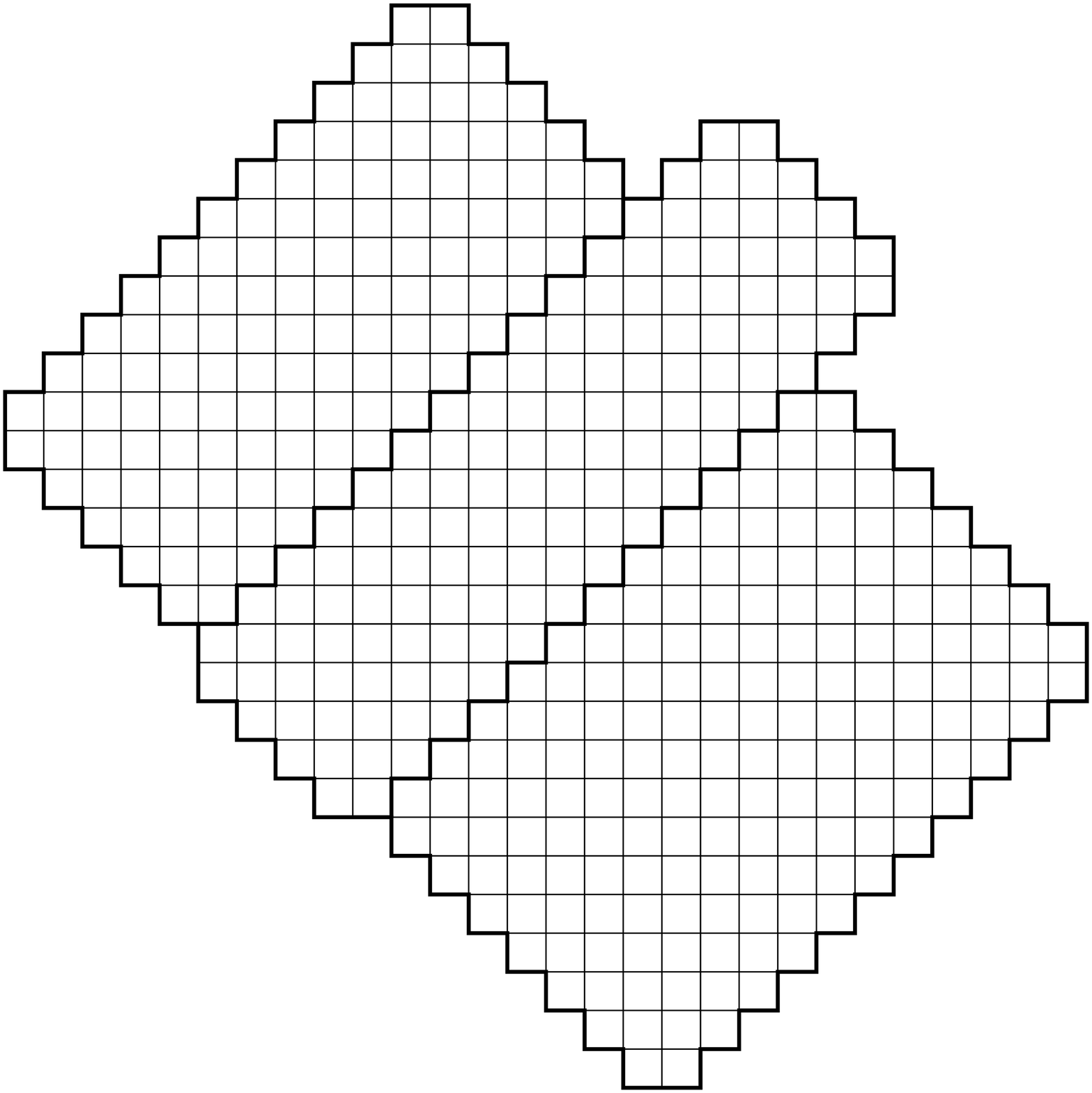}}
\hfill
{\includegraphics[width=0.53\textwidth]{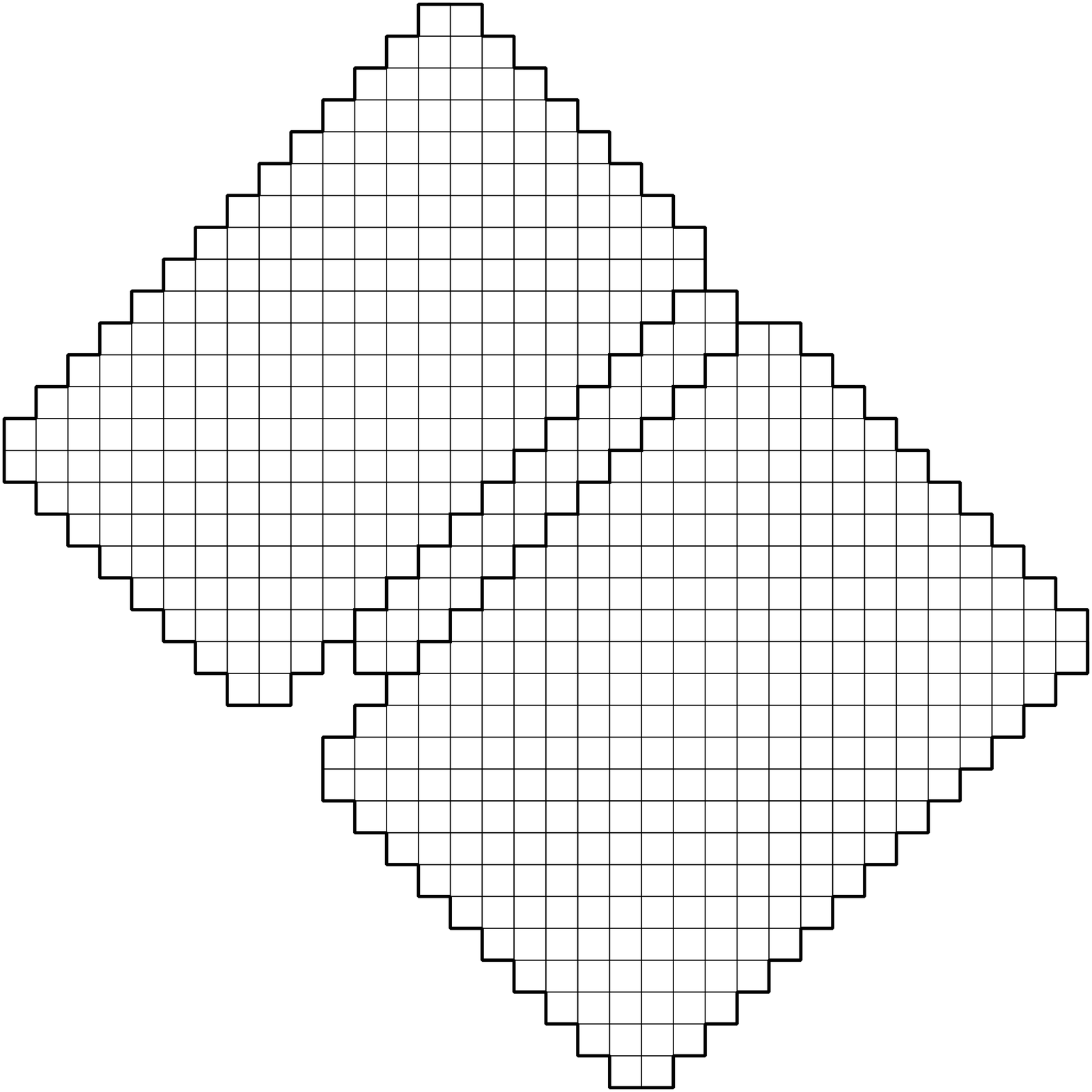}}
\hfill
}
\vskip-0.1in
\caption{The regions in Figure 4, viewed as unions of three Aztec rectangles.}
\vskip-0.1in
\label{fbd}
\end{figure}

To define $C_{m,n}^{a,b,c,d}$ when some of the parameters $a$, $b$, $c$ or $d$ are negative, let us note first what it means for any of them to equal zero: It means that the corresponding pier just sticks out by a chain of diagonally adjacent unit squares (an example with $c=0$ is illustrated in the picture on the left in Figure \ref{fbc}). Clearly, a pier of length $k>0$ can be obtained from a pier of length~0 by extending it out successively $k$ times, by inclusion of a suitable zigzag strip of unit width (this is illustrated for the northeastern pier in the picture on the left in Figure \ref{fbc}). Define then a pier of length $k<0$ by successively {\it removing} such a zigzag strip $|k|$ times from a pier of length 0. The picture on the right in Figure \ref{fbc} shows $C_{14,1}^{8,-1,10,-3}$.


Color the unit squares of the square lattice black and white in a chessboard fashion. A domino tiling of a lattice region $R$ is a covering of $R$ by horizontal or vertical dominos (unions of unit squares that share an edge) with no gaps or overlaps. Since each domino covers one black and one white unit square, the existence of a domino tiling of $R$ implies that $R$ has the same number of black and white unit squares. We call such a region balanced. Since we are interested in enumerating the domino tilings of the lattice regions we consider, unless specifically stated otherwise, throughout this paper we will assume that they are balanced.

It is not hard to show that the cruciform region $C_{m,n}^{a,b,c,d}$ is balanced if and only if
\begin{equation}
a+b+c+d=m+n-1.
\label{eba}
\end{equation}
Furthermore, it turns out\footnote{ This follows for instance from the proof of the graph splitting lemma \cite[Lemma 2.1]{CL}, as if $a>m$, the difference between the number of black and white unit squares of the pier of length $a$ (in a chessboard coloring) cannot be compensated by any choice of the dominos connecting it to the rest of the cruciform region; see also the paragraph just before the statement of Corollary \ref{tcb}.} that $C_{m,n}^{a,b,c,d}$ has no domino tilings if $a$ or $c$ is larger than $m$, or if $b$ or $d$ is larger than $n$. It follows that for a tileable cruciform region we cannot have negative values both in $\{a,c\}$ and in $\{b,d\}$. Indeed, that would imply $(a+c)+(b+d)\leq(m-1)+(n-1)<m+n-1$. Therefore, we may assume without loss of generality that $a,c\geq0$.


This allows an alternative description of cruciform regions, which is in some sense more direct. Start with an Aztec rectangle region, and place above and below it two aligned Aztec rectangles, so that the latter fit with no gaps next to the former (see Figure \ref{fbd}). The resulting region is then a cruciform region, and with this construction it is apparent that $b$ and $d$ can have negative values.

The main results of this paper are the following. Given a lattice region $R$, we denote by $\M(R)$ the number of domino tilings of $R$. The hyperfactorial function $\h(n)$ is defined by
\begin{equation}
\h(n):=0!\,1!\cdots (n-1)!
\label{ebb}
\end{equation}
\begin{theo}
\label{tba}
Let $C_{m,n}^{a,b,c,d}$ be a tileable cruciform region. Then
\begin{align}
\nonumber
\\[0pt]
\M(C_{m,n}^{a,b,c,d})&=2^{\left\{\frac14m(3m+1)+\frac14n(3n+1)+\frac12(a+c)(b+d)-\frac14(m-n)(a-b+c-d)\right\}}
\nonumber
\\[5pt]
&\ \ \ \ 
\times
\frac{\h(m+n+1)^2\h(m-a)\h(n-b)\h(m-c)\h(n-d)}{\h(n+a+1)\h(m+b+1)\h(n+c+1)\h(m+d+1)}.
\label{ebc}
\end{align}

\end{theo}

Our second main result concerns the family of regions $E_n^{a,b}$ (with $n,a,b\geq0$), called elbow regions, defined as follows. Consider the cruciform region $C_{n,n}^{a,b,b,a-1}$, and denote by $L$ its horizontal row of unit squares connecting the western and eastern outside corners. Then the elbow region $E_n^{a,b}$ is defined to consist of the portion of $C_{n,n}^{a,b,b,a-1}$ above $L$ (see the picture on the left in Figure \ref{fda} for an example). A straightforward analysis shows that $E_n^{a,b}$ is balanced if and only if $a+b=n$.

\begin{theo}
\label{tbb}
Let $E_{n}^{a,b}$ be a tileable elbow region. Then
\begin{align}
\M(E_{n}^{a,b})&=2^{n(n+1)/2}\,n!\,\frac{\h(2n+1)\h(a)\h(b)}{\h(n+a+1)\h(n+b+1)}.
\label{ebd}
\end{align}

\end{theo}

As we will see, this readily implies the following corollary.

\begin{cor}
\label{tbc}
The number of domino tilings of the region ${\mathcal T}_n$ satisfies
\begin{align}
\M({\mathcal T}_n)\,|\,\M(E_{2n-1}^{n-1,n})=2^{n(2n-1)}\,\frac{(n-1)!\,(2n-1)!}{(3n-1)!}\frac{0!\,1!\cdots(4n-2)!}{[(n-1)!\,n!\cdots(3n-2)!]^2}.
\label{ebe}
\end{align}

\end{cor}

This can be regarded as partial progress towards proving Di Francesco's conjecture that $\M({\mathcal T}_n)$ is given by formula \eqref{eaa}.

The rest of this paper is organized as follows. In Section 3 we prove Theorem \ref{tba}. Our proof is based on the complementation theorem of \cite{CT} and a result of Krattenthaler \cite{Kratt} which gives a formula for the number of perfect matchings of certain ``doubly intruded'' Aztec rectangle graphs. In Section 3 we use the factorization theorem of \cite{FT} and Theorem \ref{tba} to deduce Theorem~\ref{tbb}, and we prove Corollary \ref{tbc}. We end the paper with some concluding remarks.

\section{Proof of Theorem \ref{tba}}

Our proof of Theorem \ref{tba} is based on the complementation theorem of \cite{CT} (see Theorem 2.1 there) and a formula due to Krattenthaler \cite[Theorem 14]{Kratt} for the number of perfect matchings of Aztec rectangle graphs with certain diagonal intrusions. 

We recall here a simpler version of the complementation theorem that will suffice for our purposes. 

A finite subgraph $G$ of the grid graph $\Z^2$ is called {\it cellular} if its set of edges can be partitioned
into 4-cycles --- the {\it cells} of the graph.

A set of contiguous cells stringed up along a diagonal is called a {\it path of cells} (or simply a {\it path}). Each of the two vertices in a path that are furthest apart from one another is called an {\it extremal vertex} of $G$.

Let $X(G)$ be the set of extremal vertices of $G$. Denote by $V(G)$ the vertex set of $G$.

Let $G$ be a cellular graph, and let $H$ be a subgraph of $G$. We say that $G$ is a {\it cellular completion} 
of $H$ if

\medskip
$(i)$ $H$ is an induced subgraph of $G$
\nopagebreak

$(ii)$ $V(G)\setminus V(H) \subseteq X(G)$.

\begin{figure}[t]
\centerline{
\hfill
{\includegraphics[width=0.20\textwidth]{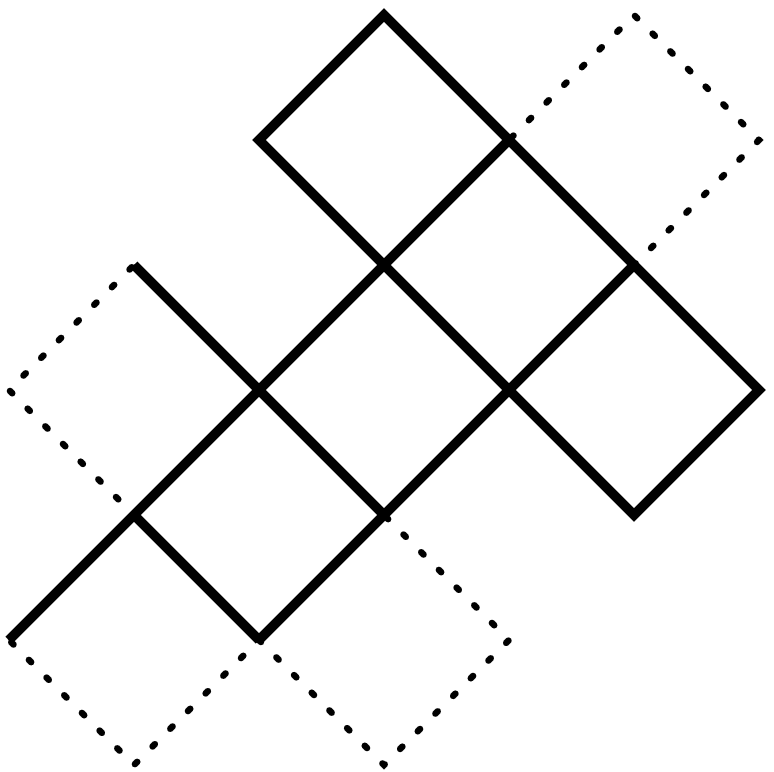}}
\hfill
{\includegraphics[width=0.20\textwidth]{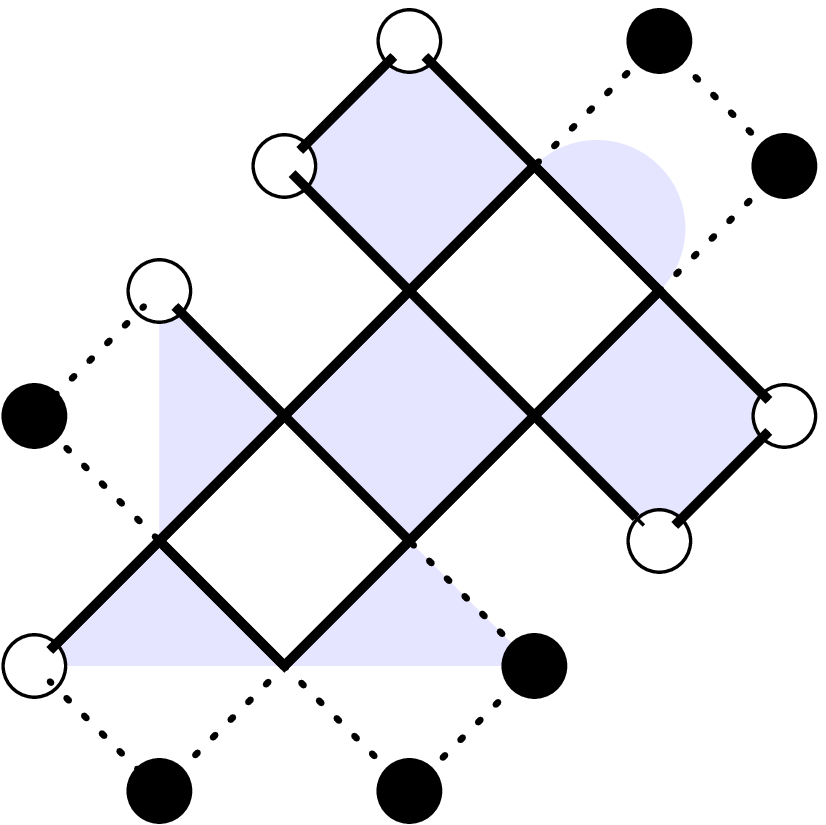}}
\hfill
{\includegraphics[width=0.20\textwidth]{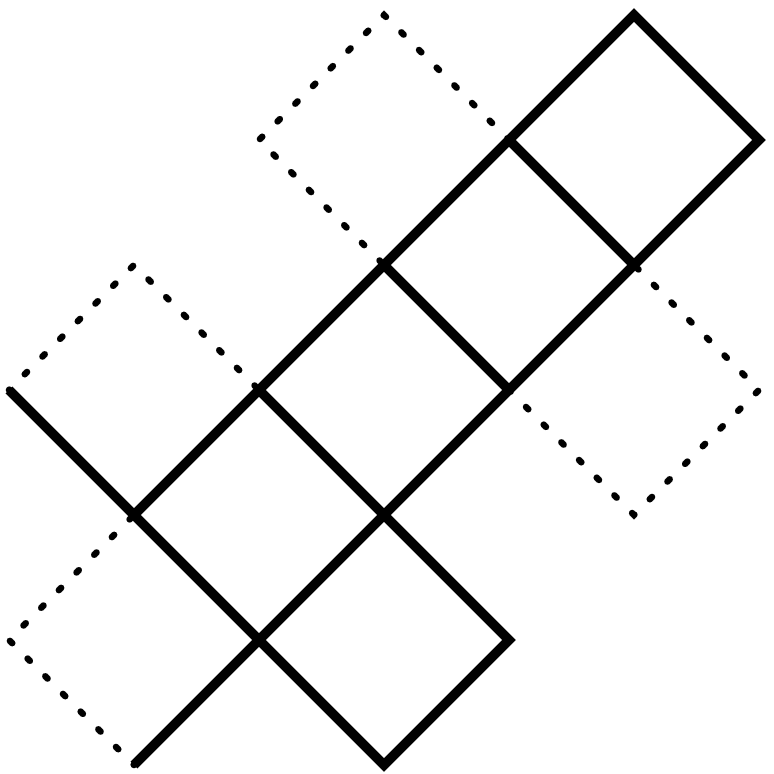}}
\hfill
}
\vskip0.1in
\caption{A graph $H$ and a cellular completion $G$ (left; the portion of $G$ that is outside $H$ is shown in dotted lines), and the complement $H'$ of $H$ with respect to $G$ (center: the white circled vertices are removed, the black circled ones are added; right: the resulting complement).}
\vskip-0.1in
\label{fcaa}
\end{figure}

\medskip
Let $G$ be a cellular completion of the graph $H$. The {\it complement
of H} (with respect to $G$) is defined to be the induced subgraph $H'$ of
$G$ whose vertex set is determined by the equation $V(H')\triangle V(H)=X(G)$,
where the triangle denotes symmetric difference of sets.
In other words, $V(H')$
is the set obtained from $V(H)$ after performing the following operation at 
each end of every path of $G$: if the corresponding extremal vertex belongs to 
$V(H)$, remove it; otherwise, include it.

An example is shown in Figure \ref{fcaa}.
The graph shown on the left in solid lines is $H$, and adding the dotted lines a cellular completion $G$ is obtained. In the center, a shading indicates the cells of $G$; the shading is ``spanned'' by the edges of $H$: a square if all four edges of the corresponding 4-cycle belong to $H$, a triangle if only two, and a semicircle if only one. The white (resp., black) circles indicate the extremal vertices of $G$ which belong (resp., do not belong) to $H$. If one discards the white vertices and adds the black ones, the resulting induced subgraph is the complement $H'$ of $H$ with respect to $G$ (this is shown on the right in the figure).


If an extremal vertex of a path $L$ of $G$ belongs to $V(H)$ then the path is
said to be {\it closed} at that end; otherwise, we say it is {\it open} at that
end. Define the {\it type} $\tau(L)$ of the path $L$ to be 1 less than the 
number of closed ends of $L$.

\begin{theo} \cite[Theorem 2.1 (Complementation Theorem)]{FT}
Let $G$ be a cellular graph with its cells partitioned into disjoint paths $L_{1},L_{2},\dotsc,L_{k}$. If $G$ is a cellular 
completion of the subgraph $H$, and $H'$ is the complement of $H$ with respect to $G$, we have
\begin{equation}
\M(H)=2^{\tau(L_{1})+\tau(L_{2})+\cdots+\tau(L_{k})}\M(H').
\label{eca}
\end{equation}
\end{theo}

\begin{figure}[t]
\centerline{
{\includegraphics[width=0.50\textwidth]{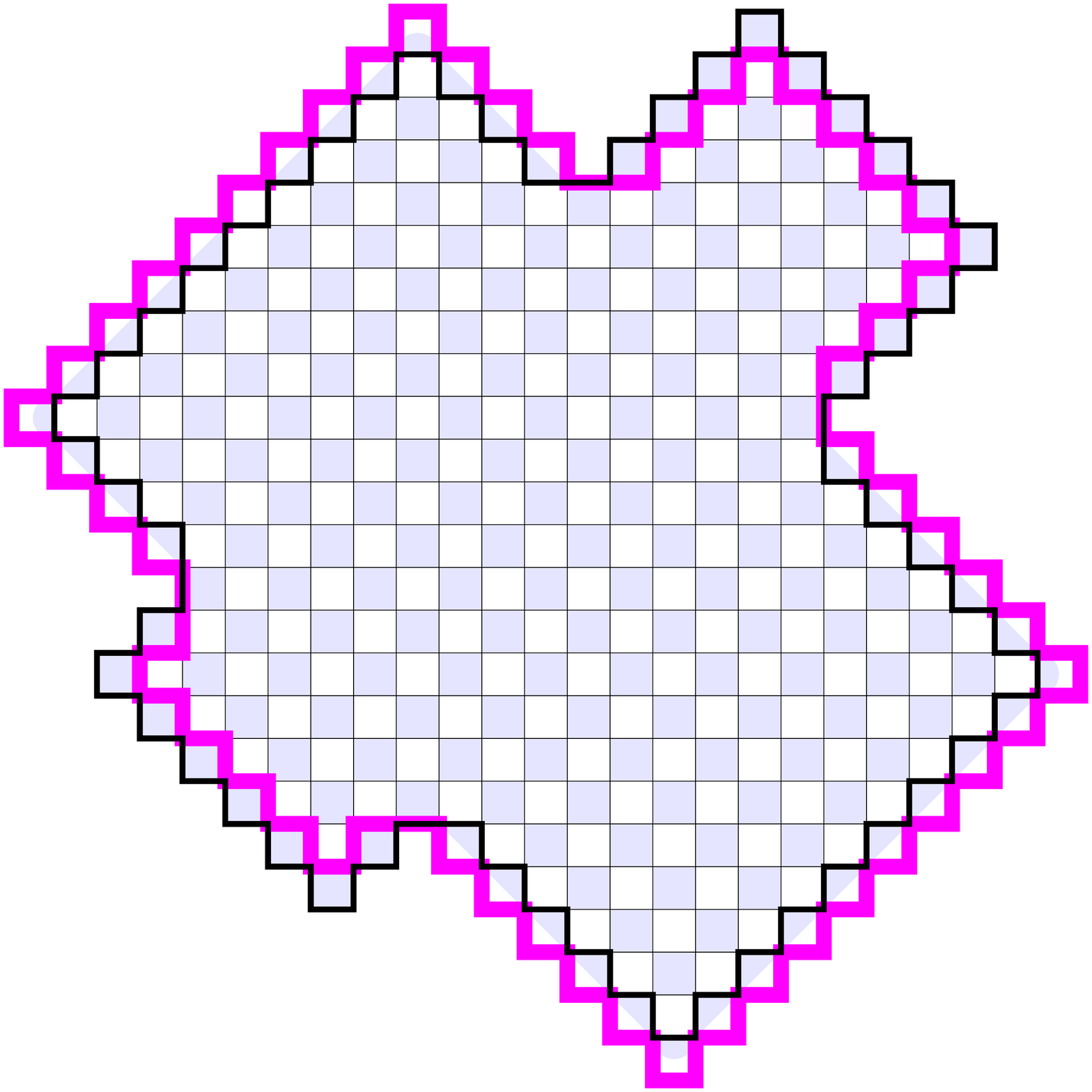}}
\hfill
{\includegraphics[width=0.55\textwidth]{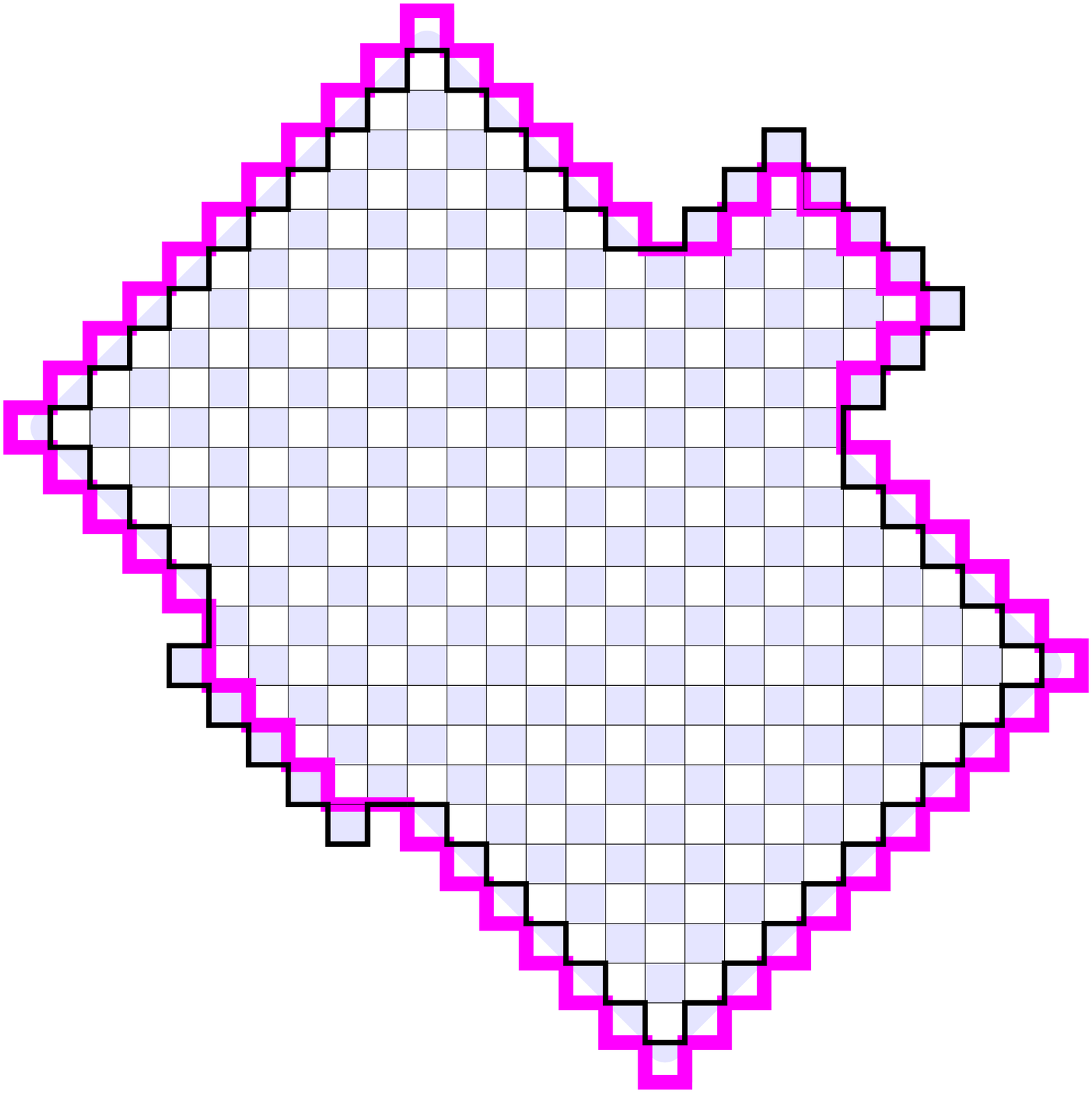}}
}
\vskip-0.1in
\caption{Applying the complementation theorem to the cruciform regions: $C_{9,6}^{3,4,5,2}$ yields $C_{9+1,6-1}^{3+1,4-1,5+1,2-1}$, i.e., $C_{10,5}^{4,3,6,1}$ (left); $C_{10,5}^{4,3,6,1}$ yields $C_{10+1,5-1}^{4+1,3-1,6+1,1-1}$, i.e., $C_{11,4}^{5,2,7,0}$ (right).}
\vskip-0.1in
\label{fca}
\end{figure}

For the situation illustrated in Figure \ref{fcaa}, the above theorem states that $\M(H)=2^0\M(H')$, as each of the three horizontal paths of cells of $G$ has precisely one closed end, hence has type~0.

As another illustration, let $H$ be the planar dual graph of the cruciform region $C_{9,6}^{3,4,5,2}$ (this is the subgraph of the grid graph induced by the lattice points on or inside the medium thickness black contour on the left in Figure \ref{fca}). The 4-cycles indicated by the light gray shading contain all its edges that are not on the boundary of the infinite face. The remaining edges are contained in ``partial cells,'' indicated by a shaded triangle (if they contain two edges of $H$) or by a shaded semicircle (if they contain a single edge of $H$).

If $G$ is the graph obtained from $H$ by including all the missing vertices and edges of these partial cells (i.e., by completing all the partial cells to 4-cycles), then $G$ is a cellular completion of $H$. The complement of $H$ with respect to $G$ is then the graph induced by the lattice points on or inside the thick contour in the figure, which happens to be\footnote{ This is actually not a coincidence, and offers in fact the key to our proof.} the planar dual graph of the cruciform region $C_{10,5}^{4,3,6,1}$. Partitioning the cells of $G$ into southwest-to-northeast going paths of cells, we see that the top 4 paths have type $\tau=-1$, the next 6 have $\tau=1$, and the last 6 have type $\tau=-1$. Therefore, by the complementation theorem we obtain
\begin{equation}
\M(C_{9,6}^{3,4,5,2})=2^{-4}\M(C_{10,5}^{4,3,6,1}).
\label{ecb}
\end{equation}

\begin{figure}[t]
\centerline{
{\includegraphics[width=0.52\textwidth]{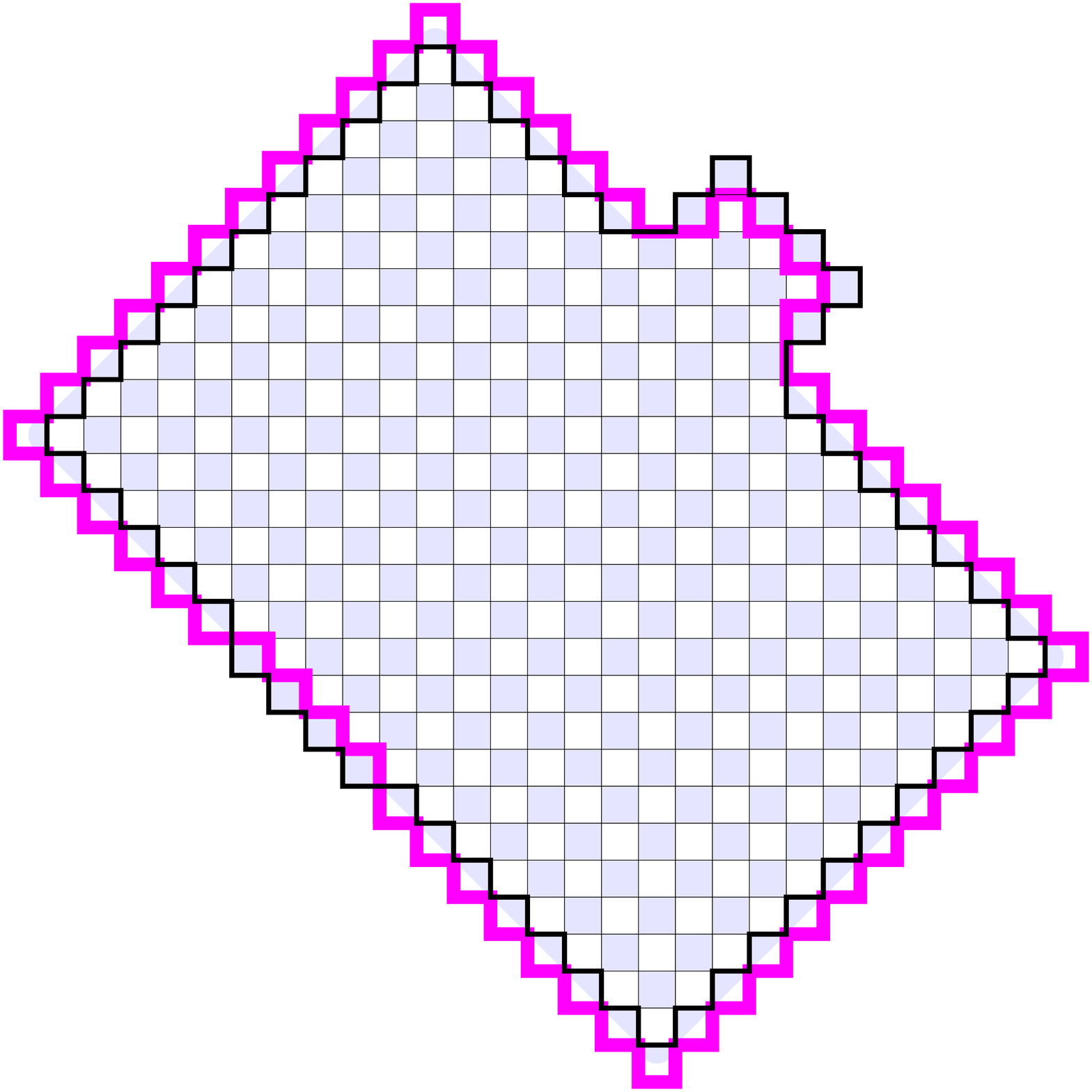}}
\hfill
{\includegraphics[width=0.55\textwidth]{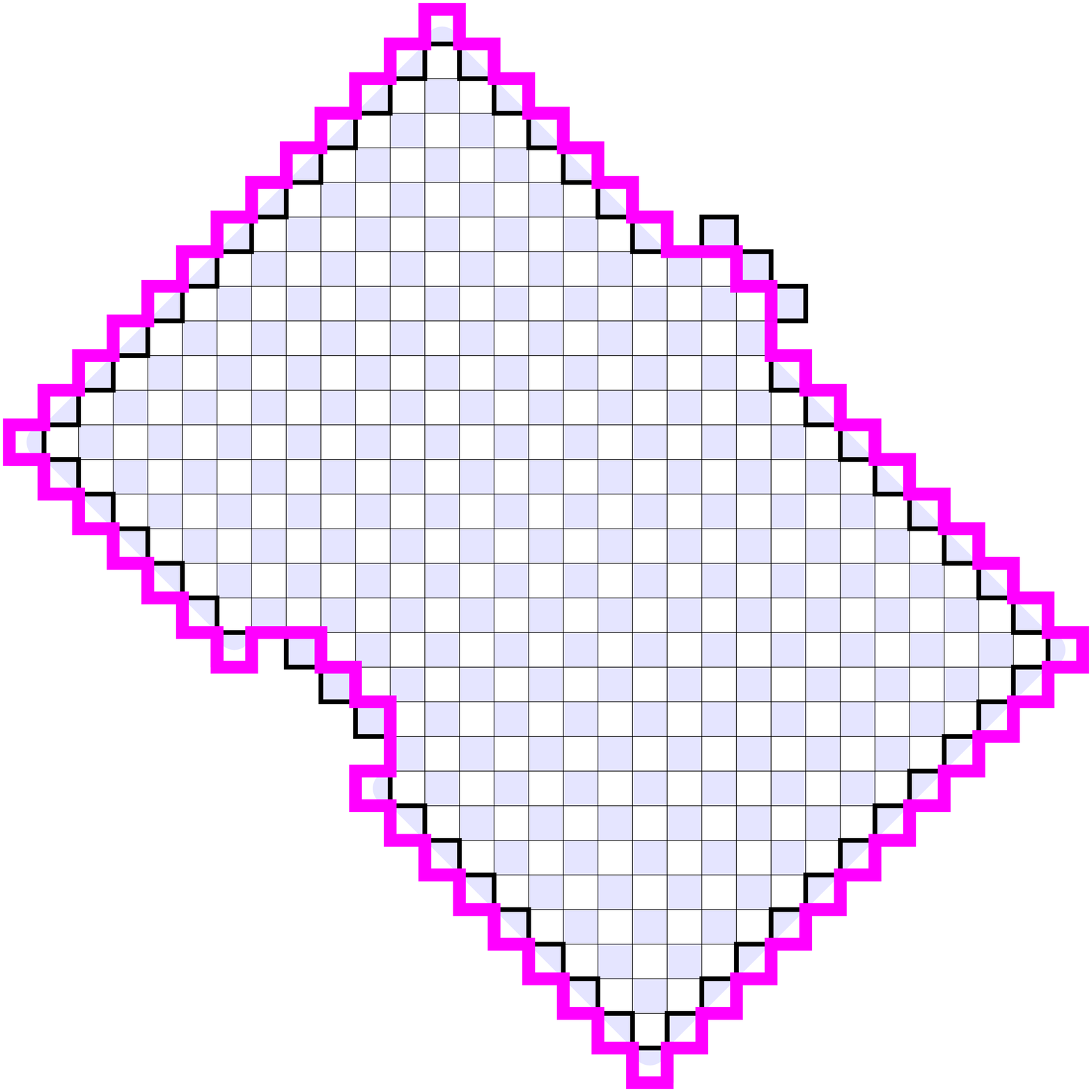}}
}
\centerline{
{\includegraphics[width=0.53\textwidth]{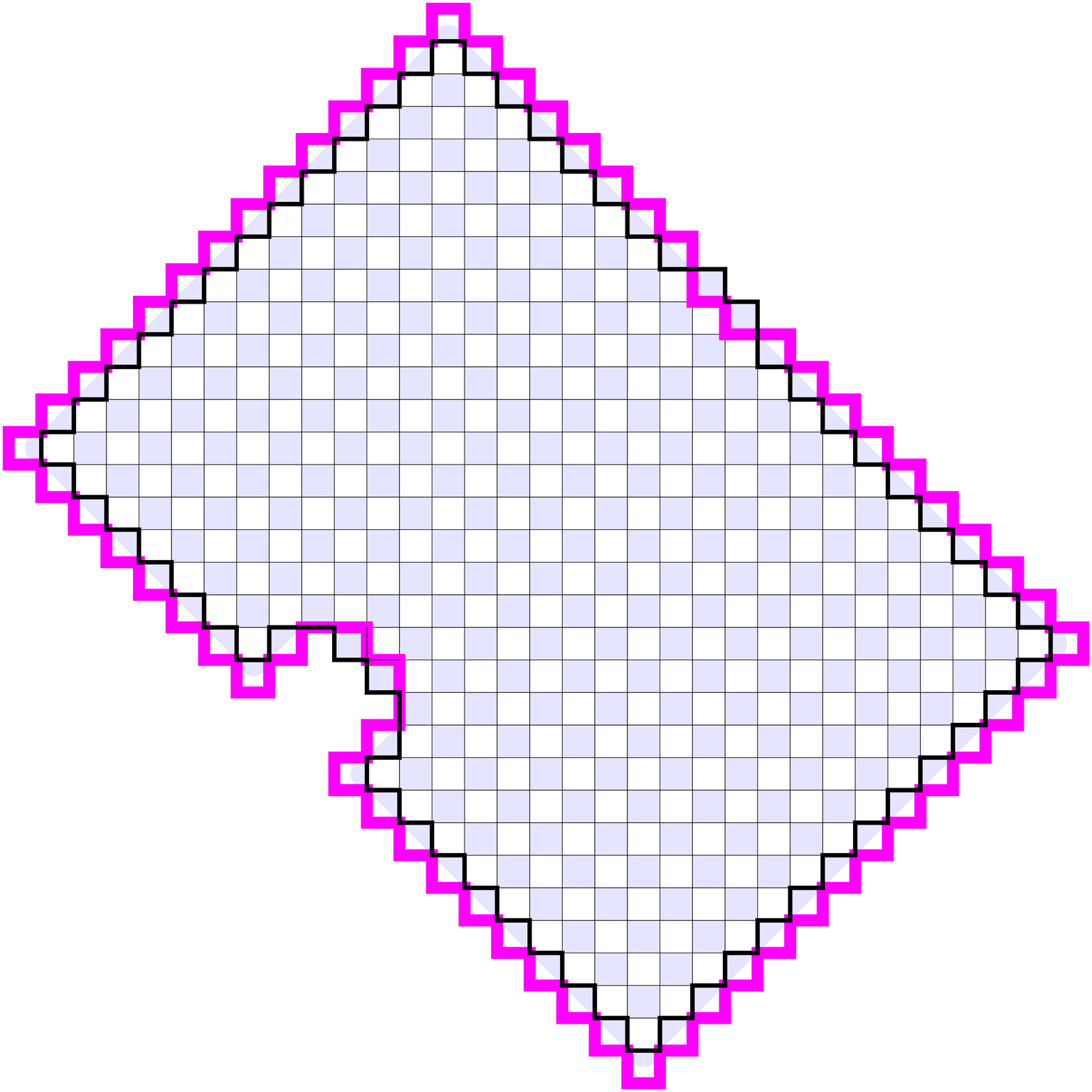}}
\hfill
{\includegraphics[width=0.53\textwidth]{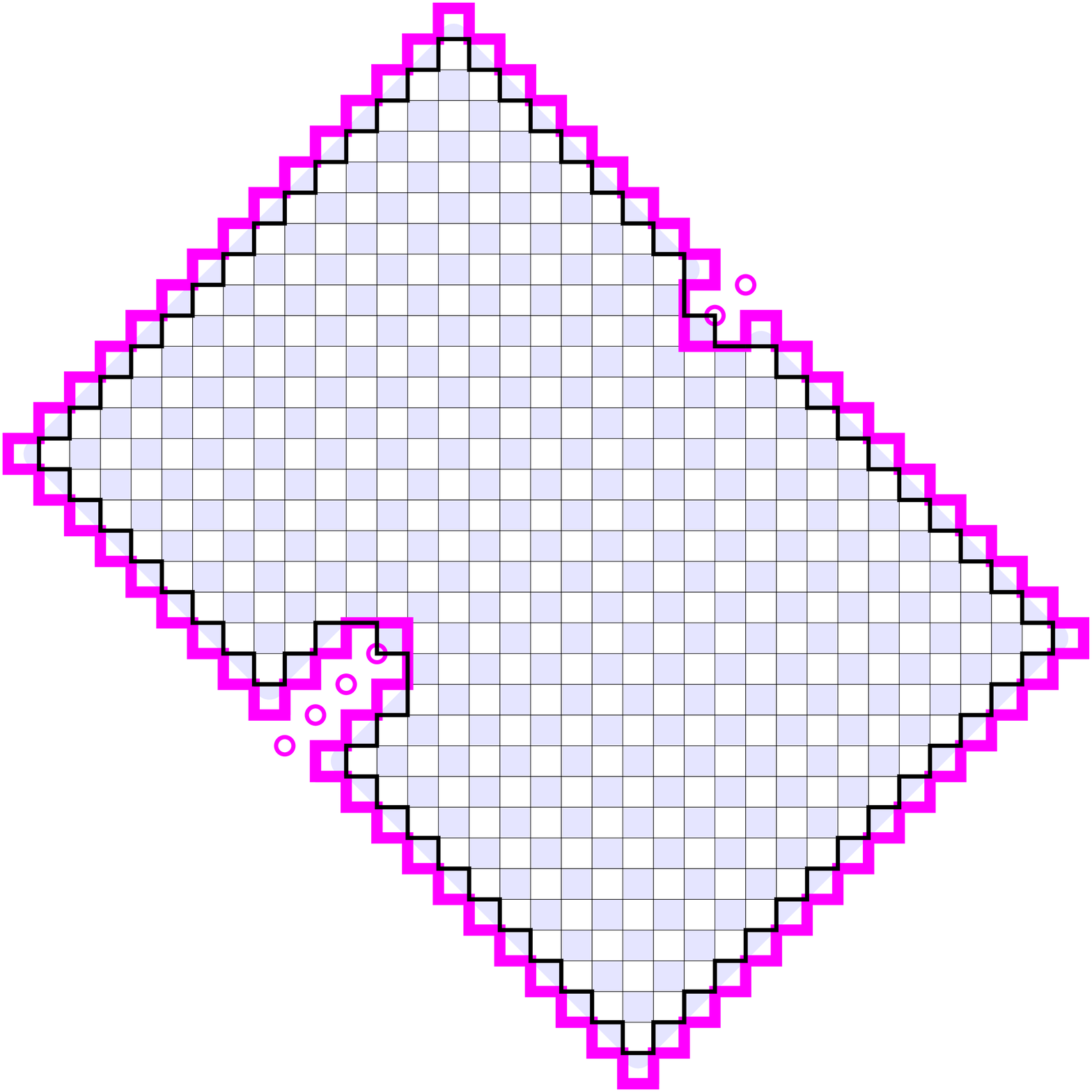}}
}
\vskip-0.1in
\caption{}
\vskip-0.1in
\label{fcb}
\end{figure}

\medskip
{\it Proof of Theorem $\ref{tba}$.} The above observation generalizes: Applying the complementation theorem this way to the planar dual graph of the cruciform region  $C_{m,n}^{a,b,c,d}$ gives rise to the dual of the cruciform region $C_{m+1,n-1}^{a+1,b-1,c+1,d-1}$. Furthermore, among the $a+n+c+2$  southwest-to-northeast going paths of cells of the former, the top $a+1$ have type $-1$, the next $n$ have type 1, and the last $c+1$ have type $-1$. Thus, the complementation theorem yields
\begin{equation}
\M(C_{m,n}^{a,b,c,d})=2^{n-a-c-2}\M(C_{m+1,n-1}^{a+1,b-1,c+1,d-1}).
\label{ecc}
\end{equation}
The same argument shows that equation \eqref{ecc} holds also if one or both of $b$ and $d$ are negative (indeed, since $a,c\geq0$ by hypothesis, the types of the southwest-to-northeast going paths of cells are then still as stated above; Figure \ref{fcb} shows three such instances). We have therefore that
\begin{equation}
\M(C_{m+i,n-i}^{a+i,b-i,c+i,d-i})=2^{n-a-c-3i-2}\M(C_{m+i+1,n-i-1}^{a+i+1,b-i-1,c+i+1,d-i-1}),
\label{ecd}
\end{equation}
for $i=0,\dotsc,n-1$ (for the dual of the cruciform region $C_{9,6}^{3,4,5,2}$, all these applications of the complementation theorem are illustrated in Figures \ref{fca} and \ref{fcb}). Combining these $n$ equations we obtain
\begin{equation}
\M(C_{m,n}^{a,b,c,d})=2^{n(n-a-c-2)-3n(n-1)/2}\M(C_{m+n,0}^{a+n,b-n,c+n,d-n}).
\label{ece}
\end{equation}

\begin{figure}[t]
\centerline{
\hfill
{\includegraphics[width=0.52\textwidth]{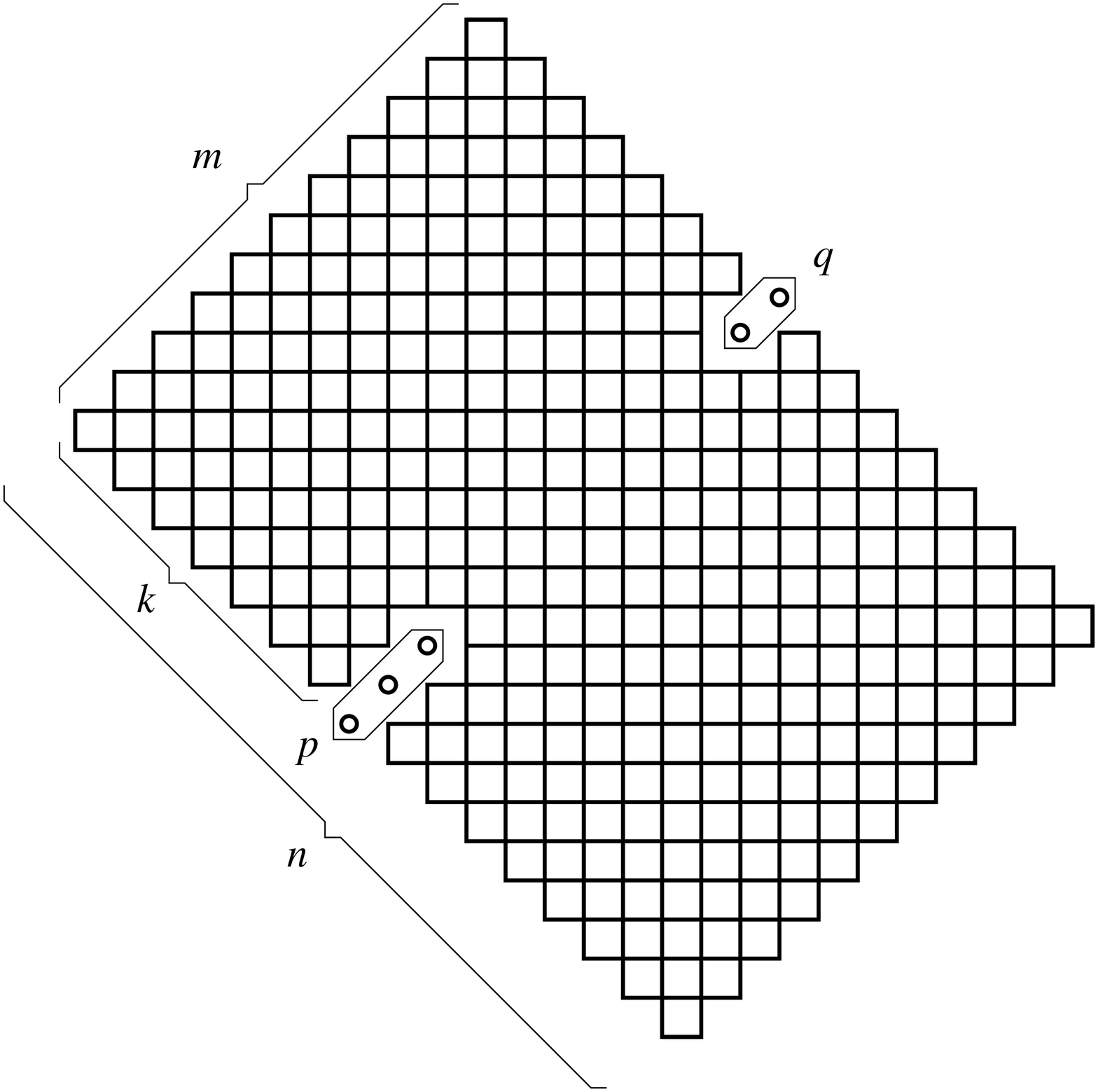}}
\hfill
{\includegraphics[width=0.52\textwidth]{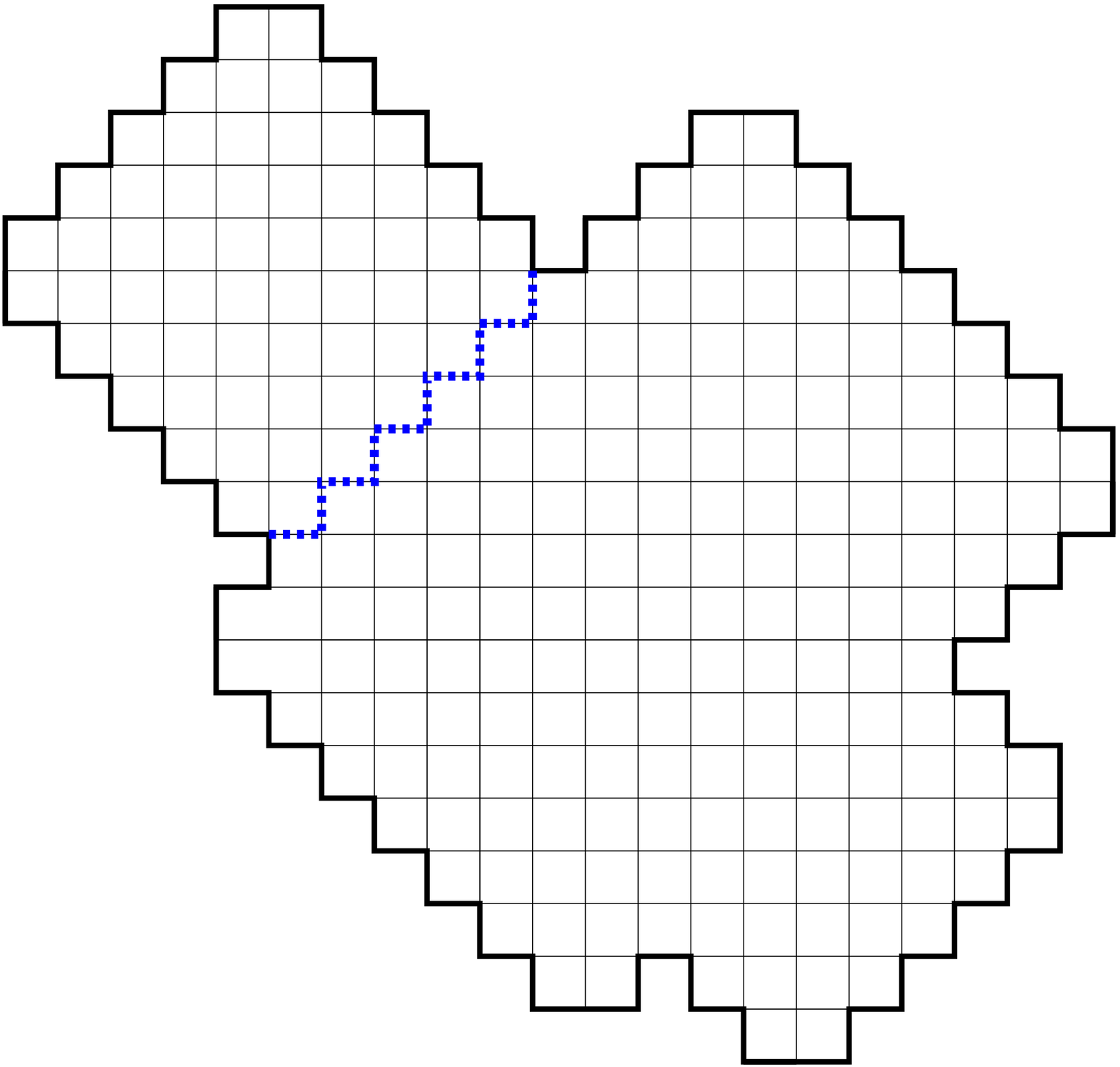}}
\hfill
}
\vskip-0.1in
\caption{{\it Left:} The doubly intruded Aztec rectangle graph $AR_{m,n}^{k}(p,q)$ for $m=11$, $n=16$, $k=7$, $p=3$, $q=2$. It is obtained from the Aztec rectangle graph $AR_{m,n}$ by removing $p$ vertices from below and $q$ vertices from above along the southwest-to-northeast diagonal whose removal
leaves $AR_{m,k}$ above it. {\it Right:} A cruciform region with $a=m$ --- $C_{5,7}^{5,3,2,1}$ --- and the $T$-region $T_{5,7}^{3,2,1}$ (the region below the dotted zigzag line).}
\vskip-0.1in
\label{fcc}
\end{figure}

It follows from the definition of the cruciform regions that the dual of $C_{m+n,0}^{a+n,b-n,c+n,d-n}$ is just the Aztec rectangle graph\footnote{ I.e., the planar dual of the Aztec rectangle region.} $AR_{m+n,2n+a+c-1}$ with some vertices along a southwest-to-northeast diagonal removed. More precisely, let $AR_{m,n}^{k}(p,q)$ be the ``doubly intruded'' Aztec rectangle graph graph shown on the left in Figure \ref{fcc}. Then the dual of $C_{m+n,0}^{a+n,b-n,c+n,d-n}$ is precisely the graph $AR_{m+n,2n+a+c+1}^{n+a}(n-d,n-b)$.

However, the number of perfect matchings of such doubly intruded Aztec rectangles follows as a special case of a result due to Krattenthaler (see \cite[Theorem 14]{Kratt}). Indeed, setting $N=m+n$, $m=n+a-1$, $d=c-a$, $C=n-d+1$ and $D=1$ in \cite[Theorem 14]{Kratt}, we obtain after some manipulation
\begin{align}
&
\M(AR_{m+n,2n+a+c+1}^{n+a}(n-d,n-b))=
2^{{2n+a+c\choose2}+(m+n+1)(m-n-a-c+1)}
\nonumber
\\[10pt]
&\ \ \ \ \ \ \ \ \ \ \ \ \ \ \ \ 
\times
\frac{\h(m+n+1)^2\h(m-a+1)\h(n-b+1)\h(m-c+1)\h(n-d+1)}{\h(n+a)\h(m+b)\h(n+c)\h(m+d)}.
\label{ecf}
\end{align}
Combining equations \eqref{ece} and \eqref{ecf}, and using $a+b+c+d=m+n-1$ (which holds since by assumption the cruciform region is tileable, hence balanced) to rewrite the resulting exponent of 2, one obtains formula \eqref{ebc}. \hfill$\square$

\medskip

The number of domino tilings of a related family of regions --- which we call $T$-regions --- can be readily deduced from Theorem \ref{tba}.

\parindent12pt

Consider the cruciform region $C_{m,n}^{m,b,c,d}$ in which the $a$-parameter is equal to $m$; an illustrative example is shown on the right in Figure \ref{fcc}. Define the $T$-region $T_{m.n}^{b,c,d}$ to be the region below the dotted zigzag line in the picture on the right in Figure \ref{fcc}. Note that by the ``graph splitting lemma'' \cite[Lemma 2.1]{CL}, the subregion above the dotted zigzag line --- which is just the Aztec diamond $AD_m$ --- must be internally tiled. It follows that
%
%
%
\begin{equation}
\M(C_{m,n}^{m,b,c,d})=\M(AD_m)\M(T_{m,n}^{b,c,d}),
\label{ecg}
\end{equation}
and, using Theorem \ref{tba} and the Aztec diamond theorem of \cite{EKLP},  we obtain the following result.    
\begin{cor}
\label{tcb}  
Let $m,n,b,c,d$ be integers with $m,n,c\geq0$ and $b+c+d=n-1$. Then the number of domino tilings of the $T$-region $T_{m.n}^{b,c,d}$ is given by
\begin{align}
\nonumber
\\[0pt]
\M(T_{m,n}^{b,c,d})&=2^{\left\{\frac14m(m-1)+\frac14n(3n+1)+\frac12(m+c)(b+d)-\frac14(m-n)(m-b+c-d)\right\}}
\nonumber
\\[5pt]
&\ \ \ \ 
\times
\frac{\h(m+n+1)\h(n-b)\h(m-c)\h(n-d)}{\h(m+b+1)\h(n+c+1)\h(m+d+1)}.
\label{ech}
\end{align}

\end{cor}




\section{Proof of Theorem \ref{tba} and Corollary \ref{tbc}}

\begin{figure}[t]
\centerline{
\hfill
{\includegraphics[width=0.385\textwidth]{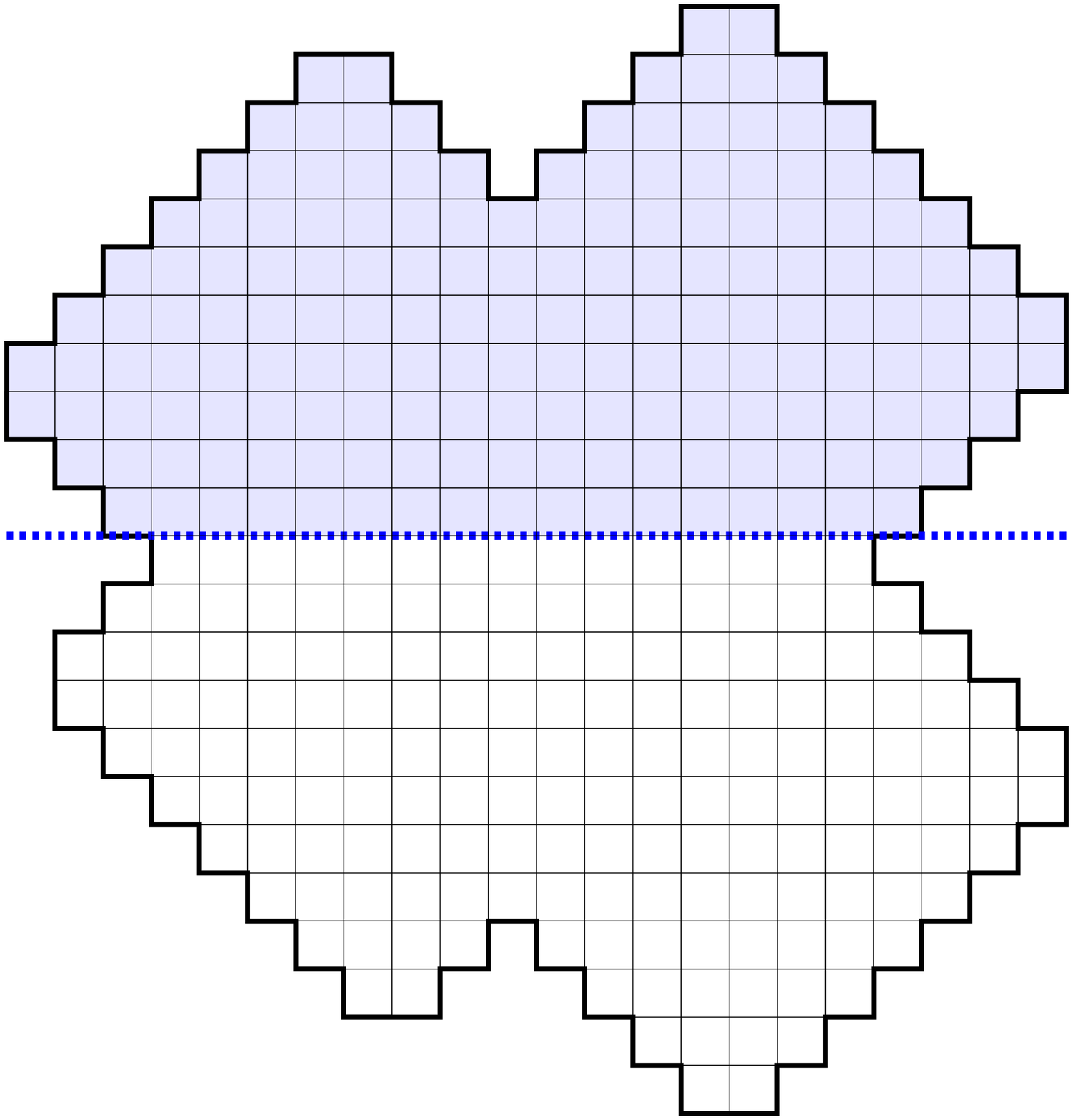}}
\hfill
{\includegraphics[width=0.45\textwidth]{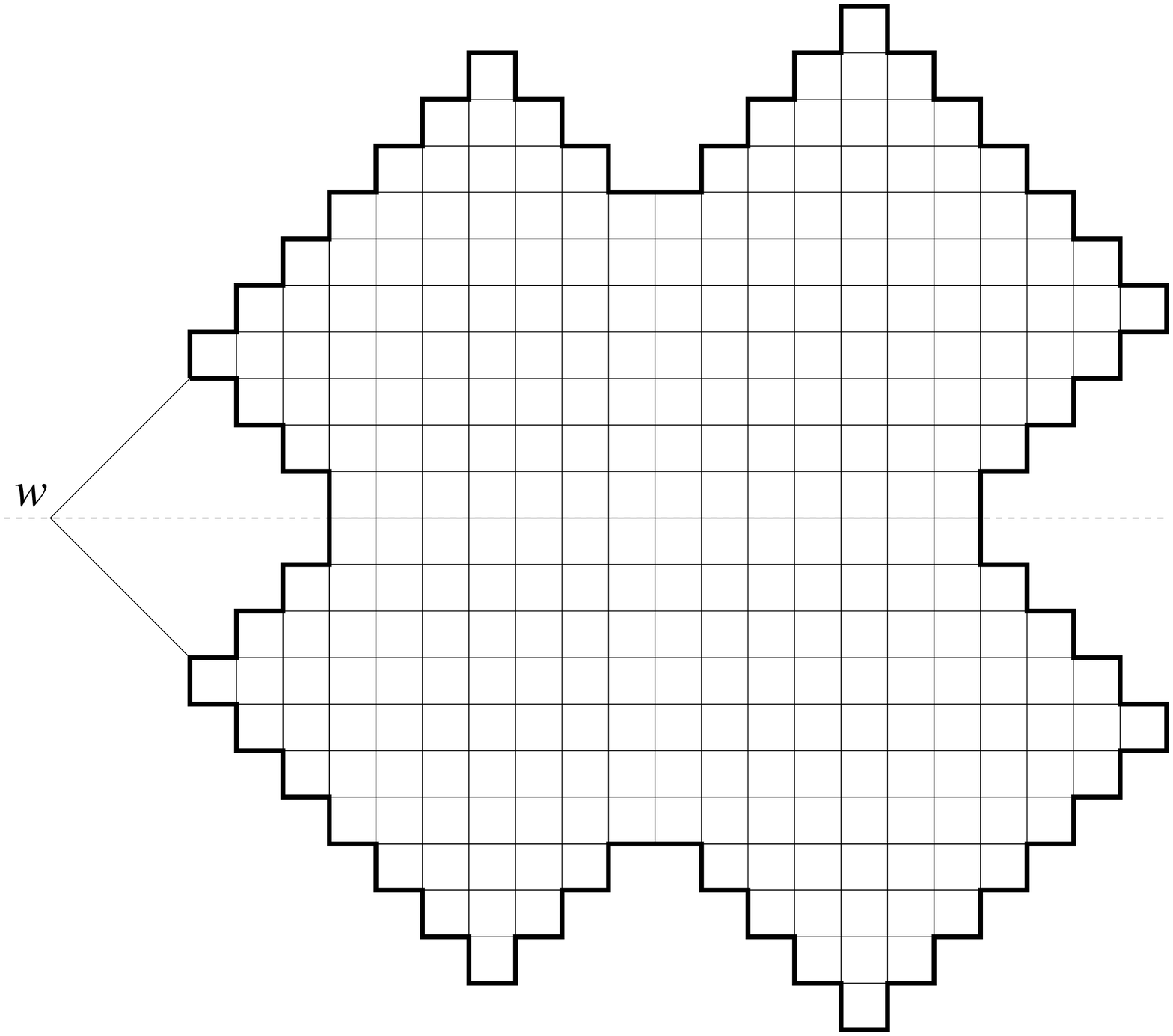}}
\hfill
}
\vskip-0.1in
\caption{{\it Left:} The elbow region $E_{n}^{a,b}$ consists of the portion of the cruciform region $C_{n,n}^{a,b,b,a-1}$ that is above its central horizontal row of unit squares (sown here is the case $n=7$, $a=3$, $b=4$). {\it Right:} The planar dual of the cruciform region $C_{7,7}^{3,4,4,3}$, with an extra vertex $w$ on the horizontal symmetry axis connected to two symmetric vertices.}
\vskip-0.1in
\label{fda}
\end{figure}

\begin{figure}[t]
\centerline{
\hfill
{\includegraphics[width=0.43\textwidth]{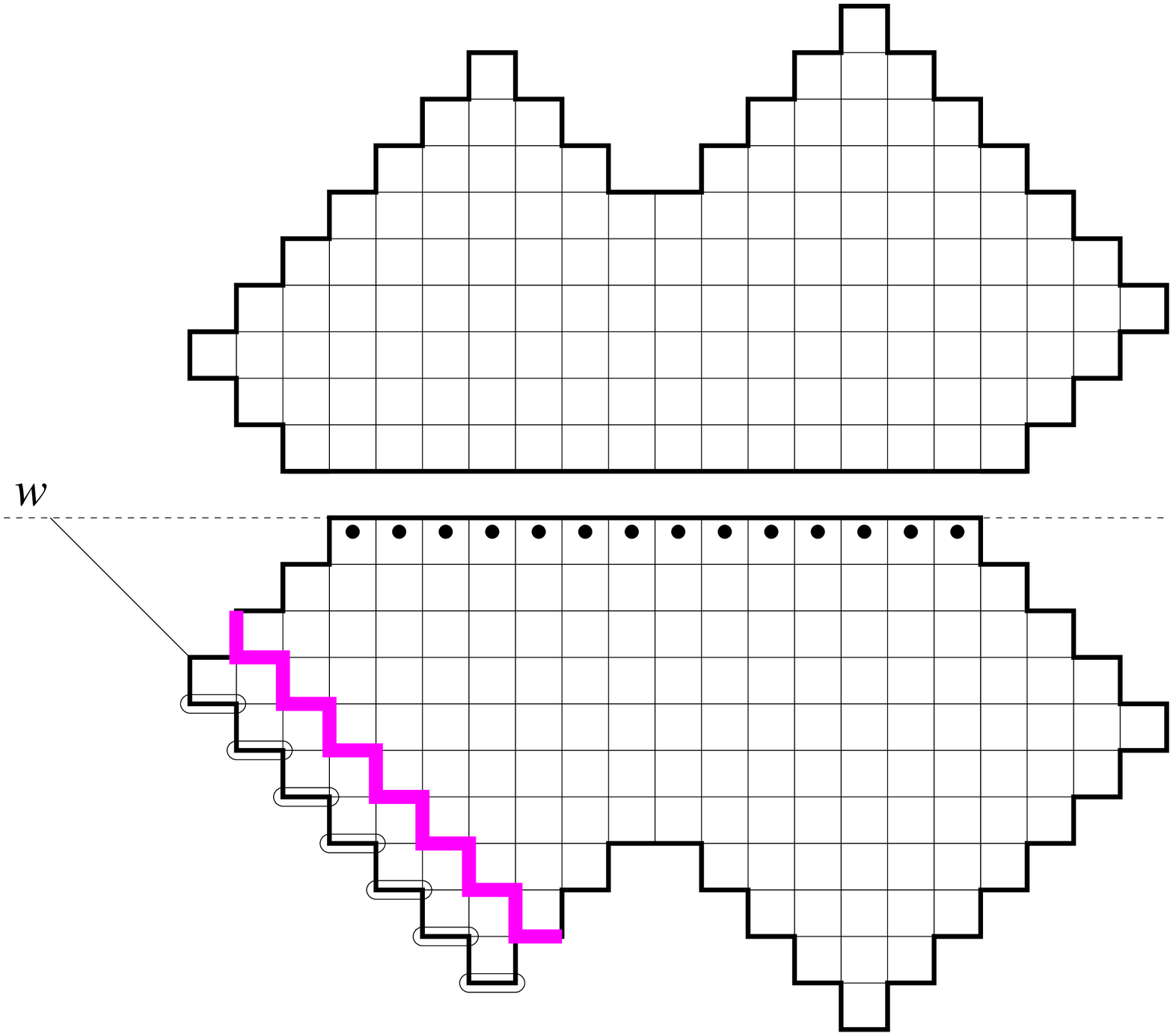}}
\hfill
{\includegraphics[width=0.47\textwidth]{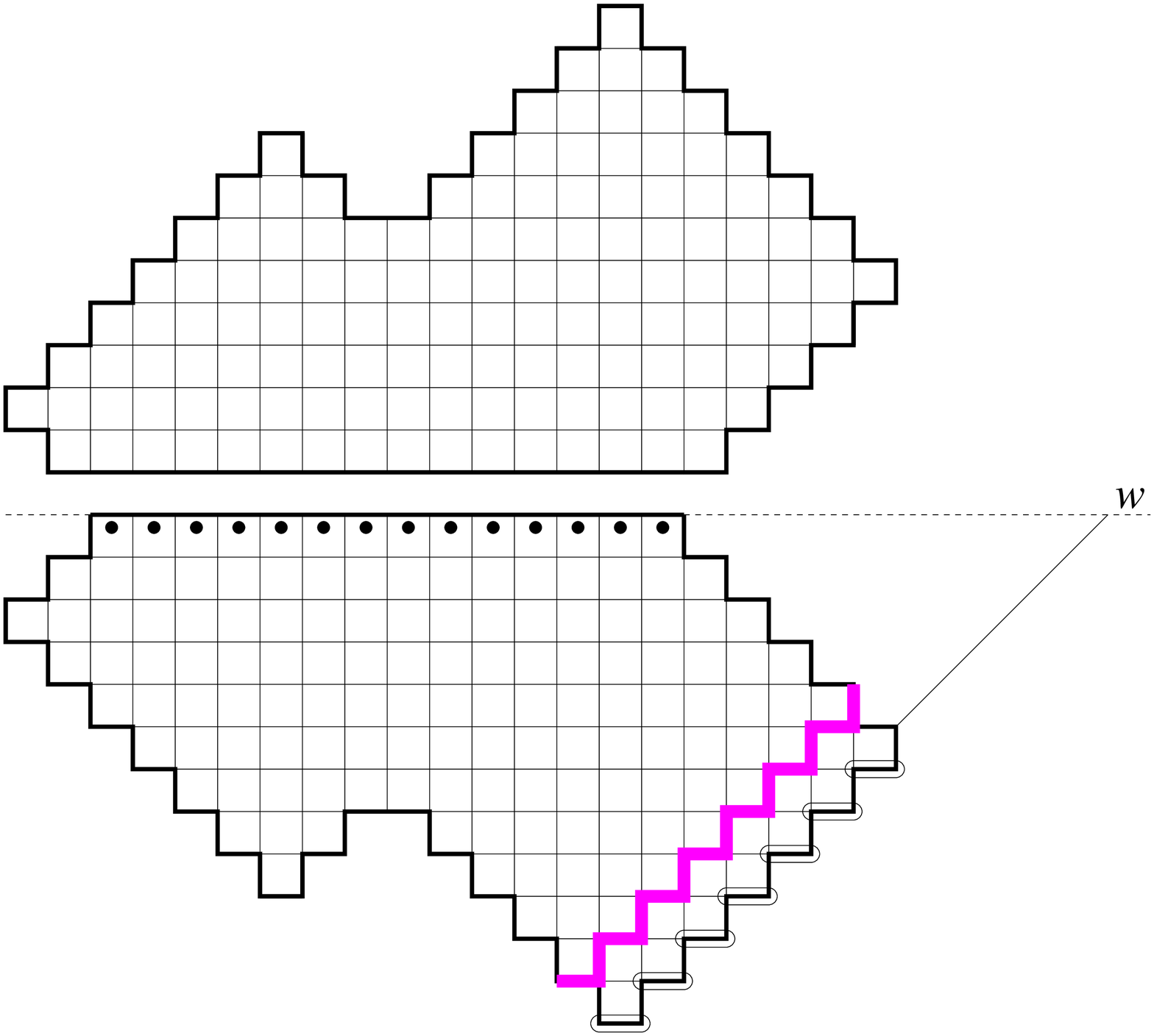}}
\hfill
}
\vskip-0.1in
\caption{Applying the factorization theorem for perfect matchings to the dual of $C_{n,n}^{a+1,b,b,a+1}$ with a vertex adjoined (left), and to the dual of $C_{n,n}^{a,b+1,b+1,a}$ with a vertex adjoined (right); in the figure, $n=7$, $a=2$, $b=4$.}
\vskip-0.1in
\label{fdb}
\end{figure}

{\it Proof of Theorem $\ref{tbb}$.}
Let $n,a,b$ be non-negative integers with $a+b+1=n$, and consider the planar dual of the cruciform region\footnote{ This cruciform graph is not balanced --- if the top leftmost vertex is white, there are one more black vertices than white ones.} $C_{n,n}^{a+1,b,b,a+1}$ (see the picture on the right in Figure \ref{fda} for an example); in order to avoid complicating the notation, we denote this graph also by $C_{n,n}^{a+1,b,b,a+1}$.

Then $C_{n,n}^{a+1,b,b,a+1}$ is a planar bipartite graph, with a horizontal axis of symmetry. It has an odd number of vertices on the symmetry axis (namely, $2n+1$), but if we extend it by including a new vertex $w$ on the symmetry axis connected to two vertices as shown in the picture on the left in Figure \ref{fdb}, the resulting graph $G$ is planar, bipartite and symmetric, with an even number of vertices on the symmetry axis. Therefore, the factorization theorem for perfect matchings \cite[Theorem 2.1]{FT} can be applied to it.

The picture on the left in  Figure \ref{fdb} shows the resulting two subgraphs. Clearly, the top one is just the dual of the elbow region $E_{n}^{a+1,b}$ (which, for simplicity of notation, we will still denote by $E_{n}^{a+1,b}$). In the bottom one --- where the dots indicate edges weighted by $1/2$ --- the single edge incident to $w$ is forced to be in every matching, which in turn, one by one, forces the same for all the circled edges. Let $F$ be the graph obtained from the bottom subgraph by removing all these forced edges and their endpoints (the southwestern boundary of $F$ is indicated by a thick gray line in the figure). Then the factorization theorem gives
\begin{equation}
\M(G)=2^{n+1}\M(E_n^{a+1,b})\M(F).
\label{eda}
\end{equation}
Clearly, in each perfect matching of $G$, the vertex $w$ must be matched either to its top or to its bottom neighbors. Since both subgraphs obtained from $G$ by removing $w$ and one of its neighbors are isomorphic to $C_{n,n}^{a+1,b,b,a}$, we have $\M(G)=2\M(C_{n,n}^{a+1,b,b,a})$, and \eqref{eda} becomes
\begin{equation}
\M(C_{n,n}^{a+1,b,b,a})=2^{n}\M(E_n^{a+1,b})\M(F).
\label{edb}
\end{equation}
In the same fashion, starting with the cruciform graph $C_{n,n}^{a,b+1,b+1,a}$, we obtain (see the picture on the right in Figure \ref{fdb})
\begin{equation}
\M(C_{n,n}^{a,b+1,b,a})=2^{n}\M(E_n^{a,b+1})\M(F),
\label{edc}
\end{equation}
where --- crucially --- the graph $F$ is precisely the same as in equation \eqref{edb}. Taking the ratio of equations \eqref{edb} and \eqref{edc} we get
\begin{equation}
\frac{\M(E_n^{a+1,b})}{\M(E_n^{a,b+1})}=
\frac{\M(C_{n,n}^{a+1,b,b,a})}{\M(C_{n,n}^{a,b+1,b,a})}.
\label{edd}
\end{equation}
Replacing $a$ by $a-1$ this becomes
\begin{equation}
\frac{\M(E_n^{a,b})}{\M(E_n^{a-1,b+1})}=
\frac{\M(C_{n,n}^{a,b,b,a-1})}{\M(C_{n,n}^{a-1,b+1,b,a-1})},
\label{ede}
\end{equation}
for all non-negative integers $a,b$ with $a+b=n$. Repeated application of \eqref{ede} yields
\begin{equation}
\frac{\M(E_n^{a,b})}{\M(E_n^{0,n})}=
\prod_{i=1}^a\frac{\M(C_{n,n}^{a-i+1,b+i-1,b+i-1,a-i})}{\M(C_{n,n}^{a-i,b+i,b+i-1,a-i})}.
\label{edf}
\end{equation}
One readily sees that, after removing the forced dominos from the elbow region $E_n^{0,n}$, the leftover region is the Aztec diamond $AD_n$, whose number of tilings is $2^{n(n+1)/2}$. Then using Theorem~\ref{tba}, equation \eqref{edf} yields formula \eqref{ebd}. \hfill$\square$

\begin{figure}[t]
\centerline{
\hfill
{\includegraphics[width=0.5\textwidth]{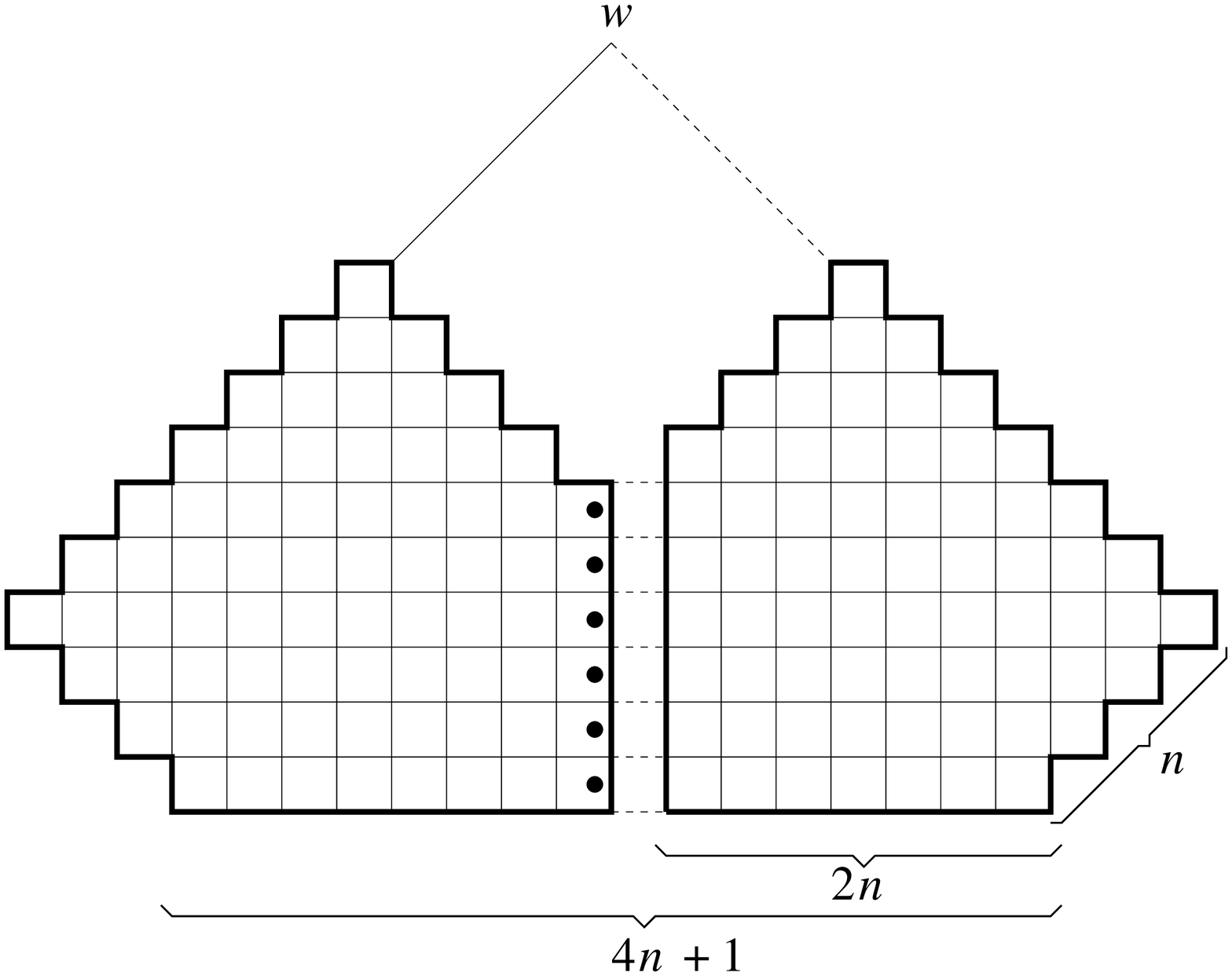}}
\hfill
{\includegraphics[width=0.275\textwidth]{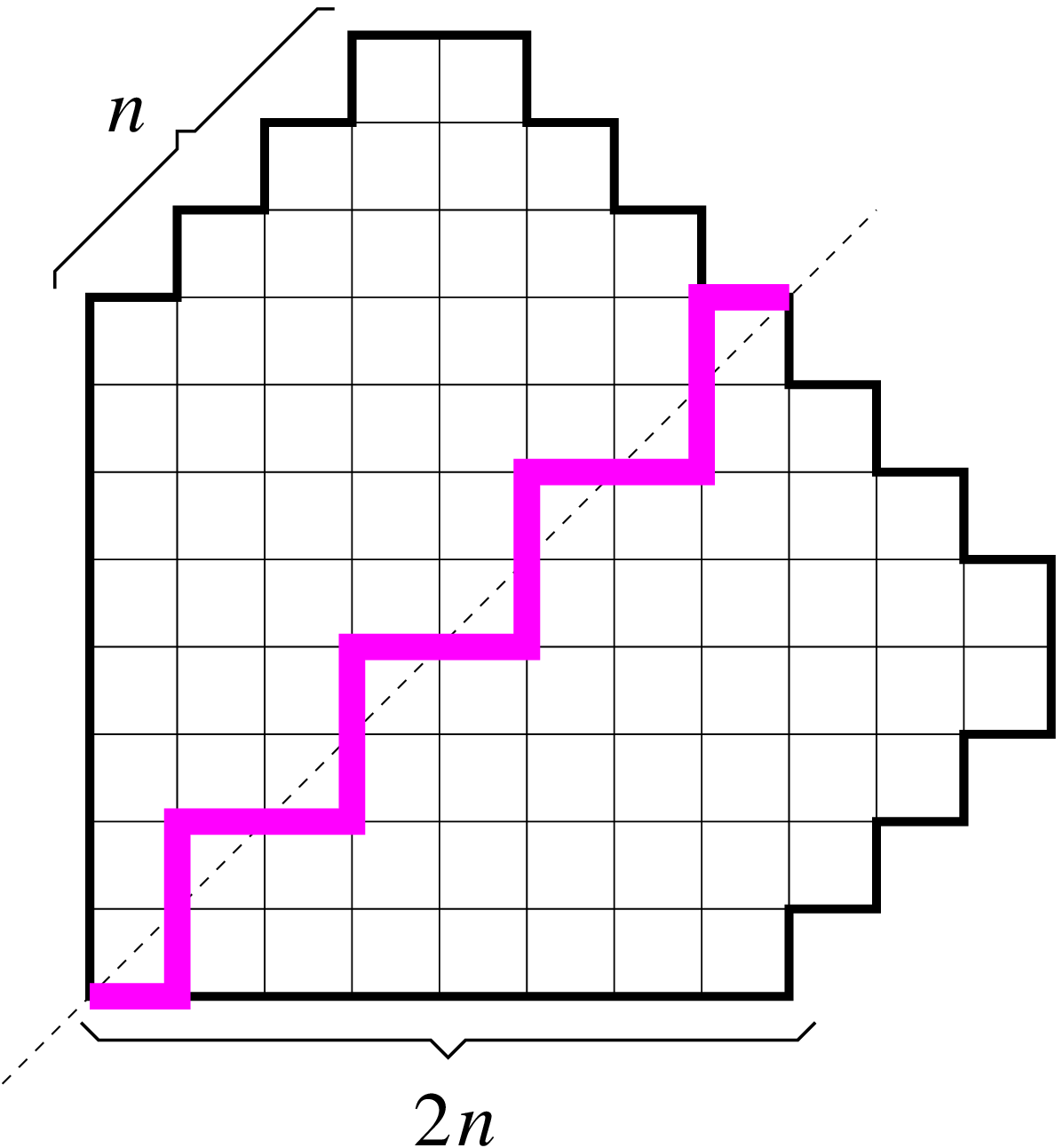}}
\hfill
}
\vskip-0.1in
\caption{Applying the factorization theorem to the dual of the elbow region $E_{2n-1}^{n,n}$ with an added vertex; the resulting subgraph on the right is $W_n$ (left). Applying the factorization theorem to the region whose dual is $W_n$ (right).}
\vskip-0.1in
\label{fdc}
\end{figure}

\medskip
{\it Proof of Corollary $\ref{tbc}$.} Consider the dual of the elbow region $E_{2n-1}^{n,n}$, and let $G$ be the graph obtained from it by adding an extra vertex $w$, joined to two original vertices as indicated on the left in Figure \ref{fdc}. Apply the factorization theorem \cite[Theorem 2.1]{FT} to $G$ (the resulting subgraphs are indicated by the thick contours in the picture on the left in Figure \ref{fdc}). Denote the resulting subgraph to the right of the symmetry axis by $W_n$, and the one to the left of the symmetry axis by $L_n$ (the $2n-2$ edges of $L_n$ marked by dots in Figure \ref{fdc} are weighted by $1/2$). Using also the fact that $\M(G)=2\M(E_{2n-1}^{n-1,n})$, we obtain
\begin{equation}
\M(E_{2n-1}^{n-1,n})=2^{n-1}\M(L_n)\M(W_n).
\label{edg}
\end{equation}
Since a perfect matching of $L_n$ can contain at most $n-1$ edges weighted by $1/2$, $2^{n-1}\M(L_n)$ is an integer. Therefore, \eqref{edg} shows that $\M(W_n)$ is a divisor of $\M(E_{2n-1}^{n-1,n})$.

Apply now the factorization theorem to the graph $W_n$. Since the symmetry axis is a diagonal of the square grid, the factorization theorem is especially simple to express in terms of domino tilings of the region whose dual is the graph $W_n$, which is represented in the picture on the right in Figure \ref{fdc} for $n=4$. Namely, the factorization theorem states in this case that $\M(W_n)$ is just the product of the number of domino tilings of the two regions obtained by cutting along the indicated zigzag path. Since the region above is just ${\mathcal T}_n$, this implies that $\M({\mathcal T}_n)$ is a divisor of $\M(W_n)$, and hence of $\M(E_{2n-1}^{n-1,n})$. \hfill$\square$

\section{Concluding remarks}

In this paper we proved a simple product formula for the number of domino tilings of certain cruciform shaped regions that generalize Aztec diamonds. They are determined by closed contours on the square grid in which all boundaries are zigzags, and one is only allowed two kinds of turns --- $90^\circ$ right turns of the kind encountered when traveling the boundary of an Aztec diamond region clockwise, and $90^\circ$ left turns of the kind shown on the right in Figure~\ref{fba} --- with the additional condition that if one encounters two consecutive turns of the same type, then the zigzag leading to the first and the zigzag leaving the second must have the same length. This description is reminiscent of the shamrock regions considered in \cite{vf}.

It would be interesting to investigate the limit shape of a random tiling of these cruciform regions. They are determined by the occupation probabilities of individual dominos, which by definition are fractions whose denominator is the number of tilings of the region, which we determined in this paper.






\end{document}